\documentclass[12pt,a4paper]{amsart}

\usepackage{euscript,amsfonts,amssymb,amsmath,amscd,amsthm,enumerate,hyperref}

\usepackage{comment}

\usepackage{mathrsfs}

\usepackage{graphicx}

\usepackage{color}

\usepackage[margin=1.1in]{geometry}

\usepackage{bbm}

\newtheorem{theorem}{Theorem}[section]
\newtheorem*{theorem*}{Theorem}
\newtheorem{lemma}[theorem]{Lemma}
\newtheorem{definition}[theorem]{Definition}
\newtheorem{proposition}[theorem]{Proposition}

\newtheorem{assumption}[theorem]{Assumption}

\theoremstyle{remark}
\newtheorem{rmk}[theorem]{Remark}

\newcommand{\E}{\mathbb E}

\newcommand{\pp}{\mathbb{P}}

\newcommand{\rr}{\mathbb{R}}

\newcommand{\nn}{\mathbb{N}}
\newcommand{\ep}{\hfill \ensuremath{\Box}}
\newcommand{\eq}{\begin{equation}}
\newcommand{\en}{\end{equation}}
\newcommand{\Mf}{M_{\mathrm{fin}}(\rr)}
\newcommand{\Mft}{M_{\mathrm{fin}}([0,t]\times\rr)}

\numberwithin{equation}{section}

\title[Fluctuations in rank-based models]{SPDE limit of the global fluctuations in rank-based models} 

\author{Praveen Kolli}
\address{Department of Mathematical Sciences, Carnegie Mellon University, PA, USA}
\email{kpc.0915@gmail.com}

\author{Mykhaylo Shkolnikov}
\address{ORFE Department, Princeton University, Princeton, NJ, USA}
\email{mshkolni@gmail.com}
\thanks{Research supported in part by NSF grant DMS-1506290.}

\begin{document}

\maketitle

\begin{abstract}
We consider systems of diffusion processes (``particles'') interacting through their ranks (also referred to as ``rank-based models'' in the mathematical finance literature). We show that, as the number of particles becomes large, the process of fluctuations of the empirical cumulative distribution functions converges to the solution of a linear parabolic SPDE with additive noise. The coefficients in the limiting SPDE are determined by the hydrodynamic limit of the particle system which, in turn, can be described by the porous medium PDE. The result opens the door to a thorough investigation of large equity markets and investment therein. In the course of the proof we also derive quantitative propagation of chaos estimates for the particle system.
\end{abstract}

\section{Introduction}

We study systems of interacting diffusion processes (``particles'') on the real line whose dynamics are given by the SDEs
\begin{equation}\label{ips_sde}
\mathrm{d}X^{(n)}_i(t)=b\big(F_{\rho^{(n)}(t)}\big(X^{(n)}_i(t)\big)\big)\,\mathrm{d}t + \sigma\big(F_{\rho^{(n)}(t)}\big(X^{(n)}_i(t)\big)\big)\,\mathrm{d}B^{(n)}_i(t),\quad i=1,\,2,\,\ldots,\,n.
\end{equation}
Here $b$, $\sigma$ are functions from $[0,1]$ to $\rr$, $(0,\infty)$, respectively, $\rho^{(n)}(t):=\frac{1}{n}\sum_{i=1}^n \delta_{X^{(n)}_i(t)}$ is the empirical measure of the particle system at time $t$, $F_{\rho^{(n)}(t)}$ is the cumulative distribution function of $\rho^{(n)}(t)$, and $B^{(n)}_1,\,B^{(n)}_2,\,\ldots,\,B^{(n)}_n$ are independent standard Brownian motions. Note that the drift and diffusion coefficients of a process $X^{(n)}_i$ take the values $b\big(\frac{k}{n}\big)$ and $\sigma\big(\frac{k}{n}\big)$ whenever the \textit{rank} (from the left) of $X^{(n)}_i(t)$ within $\big(X^{(n)}_1(t),X^{(n)}_2(t),\ldots,X^{(n)}_n(t)\big)$ is $k$. This allows to identify \eqref{ips_sde} with the so-called rank-based models of stochastic portfolio theory introduced by \textsc{Fernholz} and \textsc{Karatzas} (see \cite[Section 13]{FK}). 

\medskip

Rank-based models have recently received much attention in pure and applied probability theory. Originally, they appeared as a special case in the context of the piecewise linear filtering problem in \cite{BP} where weak uniqueness for \eqref{ips_sde} is established (weak existence being a consequence of the general result in \cite[Exercise 12.4.3]{SV}). The recent renewed interest in rank-based models stems from the fact that they are the first ones to capture the shape and stability of the capital distribution among companies in the U.S. We refer to \cite[Figure 5.1]{Fe} for a plot of the U.S. capital distribution curves over seventy years and to \cite{CP} and \cite{IPS} for the mathematical results on their shape and stability in the setting of rank-based models. In this context, one is particularly interested in the large $n$ behavior of the system \eqref{ips_sde} which describes the evolution of the capital distribution when one takes thousands of companies into account. The stocks of the latter comprise typical portfolios of institutional investors and the change in the capital distribution is central to their investment decisions.

\medskip

We point out that \eqref{ips_sde} falls into the general framework of particle systems interacting through their mean field whose analysis originates with the seminal work \cite{Mc} of \textsc{McKean}. In the case of diffusion processes the general results on the subject can be summarized as follows. A law of large numbers as $n\to\infty$ (``hydrodynamic limit'') has been obtained assuming the joint continuity of the drift and diffusion coefficients with respect to the current location of the particle and the empirical measure by \textsc{G\"artner} in \cite{Ga} (see also \cite{Le}, \cite{Oe1} for previous results under more restrictive assumptions). Gaussian fluctuations around the hydrodynamic limit have been established for drift coefficients $\int_\rr b(X^{(n)}_i(t),y)\,\rho^{(n)}(t)(\mathrm{d}y)$ with a twice continuously differentiable function $b$ and constant diffusion coefficients by \textsc{Tanaka} in \cite{Ta} (see also \cite{Oe2} for the case of $\rr^d$ and drift coefficients of gradient type). Concurrently, \textsc{Sznitman} \cite{Sz1} proved the Gaussian nature of the fluctuations in the absense of drift and with diffusion coefficients $\int_\rr \sigma(X^{(n)}_i(t),y)\,\rho^{(n)}(t)(\mathrm{d}y)$ with a twice continuously differentiable function $\sigma$. Finally, large deviations around the hydrodynamic limit have been studied by \textsc{Dawson} and \textsc{G\"artner} \cite{DG} in the case of a jointly continuous drift coefficient and a continuous diffusion coefficient depending only on the current location of the particle. 

\medskip

None of the described results can be applied to the system \eqref{ips_sde} due to the discontinuity of both the drift and the diffusion coefficients. Nonetheless, the special structure of the coefficients in \eqref{ips_sde} made it possible to derive the hydrodynamic limit of that system (see \cite[Proposition 2.1]{JR} and also \cite[Corollary 1.6]{DSVZ}, \cite[Theorem 1.2]{S}). More specifically, let $M_1(\rr)$ be the space of probability measures on $\rr$ equipped with the topology of weak convergence and $C([0,\infty),M_1(\rr))$ be the space of continuous functions from $[0,\infty)$ to $M_1(\rr)$ endowed with the topology of locally uniform convergence. Given that the initial positions $X^{(n)}_1(0),X^{(n)}_2(0),\ldots,X^{(n)}_n(0)$ are i.i.d. according to a probability measure $\lambda$ with a finite first moment and that $b$ and $\sigma$ in \eqref{ips_sde} are continuous, the functions $t\mapsto\rho^{(n)}(t)$, $n\in\nn$ converge in probability in $C([0,\infty),M_1(\rr))$ to a deterministic limit $t\mapsto\rho(t)$. Moreover, the associated cumulative distribution functions $R(t,\cdot):=F_{\rho(t)}(\cdot)$, $t\ge0$ form the generalized solution to the Cauchy problem for the \textit{porous medium equation}:   
\begin{equation}\label{PME}
R_t=-B(R)_x+\Sigma(R)_{xx},\quad R(0,\cdot)=F_\lambda(\cdot),
\end{equation}
where $B(r):=\int_0^r b(a)\,\mathrm{d}a$ and $\Sigma(r):=\int_0^r \frac{1}{2}\,\sigma(a)^2\,\mathrm{d}a$ (\cite[Definition 3]{Gi} of a generalized solution to \eqref{PME} is briefly recalled in Definition \ref{Gi_def} below). In fact, under additional moment and regularity assumptions it has been shown in \cite[Theorem 1.4]{DSVZ} that the sequence $t\mapsto\rho_n(t)$, $n\in\nn$ satisfies a large deviation principle in $C([0,\infty),M_1(\rr))$. 

\medskip

In this paper we are concerned with the \textit{fluctuations} of the particle system \eqref{ips_sde}. To this end, we introduce the space $\Mf$ of finite signed measures on $\rr$, viewed as the dual of $C_0(\rr)$ and endowed with the associated weak-$*$ topology. Similarly, we define the spaces $\Mft$ for $t>0$ and equip each of them with the respective weak-$*$ topology. The fluctuations of the particle system \eqref{ips_sde} are studied via the $\Mf$-valued processes 
\begin{equation}\label{def:Gn}
t\mapsto G_n(t)(\mathrm{d}x):=\sqrt{n}\,(F_{\rho^{(n)}(t)}(x)-R(t,x))\,\mathrm{d}x,\quad n\in\nn
\end{equation} 
indexed by $t\in[0,\infty)$, as well as the processes
\begin{equation}\label{def:Hn}
t\mapsto H_n(t)(\mathrm{d}s,\mathrm{d}x):=\sqrt{n}\,(F_{\rho^{(n)}(s)}(x)-R(s,x))\,\mathrm{d}x\,\mathrm{d}s,\quad n\in\nn
\end{equation}
taking values in $\Mft$, $t>0$, respectively. Note that the measures $G_n(t)$, $t\ge0$ belong to $\Mf$ and the measures $H_n(t)$, $t>0$ are elements of $\Mft$, $t>0$ as soon as the first moments of the probability measures $\rho(t)$, $t\ge0$ are finite and uniformly bounded on compact intervals of $t$'s. This turns out to be the case under the following assumption (see the estimate \eqref{moment_est} below). 

\begin{assumption}\label{main_asmp}
\begin{enumerate}[(a)]
\item There exist $\eta>0$ and $\lambda\in M_1(\rr)$ such that $\lambda$ has a bounded density and finite moments up to order $(2+\eta)$ and the initial positions $X^{(n)}_1(0),X^{(n)}_2(0),\ldots,X^{(n)}_n(0)$ are i.i.d. according to $\lambda$ for all $n\in\nn$. 
\item The functions $b$ and $\sigma$ in \eqref{ips_sde} are differentiable with locally H\"older continuous derivatives.
\end{enumerate}
\end{assumption}

\smallskip

Before proceeding it is worth it to point out that the processes of \eqref{def:Gn} and \eqref{def:Hn} provide access to observables of the forms 
\begin{eqnarray}
&&\;\;\int_\rr \gamma(x)\,G_n(t)(\mathrm{d}x)=\sqrt{n} \int_\rr \bigg(\!\int_0^x\! \gamma(y)\,\mathrm{d}y\!\bigg) (\rho^{(n)}(t)(\mathrm{d}x)\!-\!\rho(t)(\mathrm{d}x)), \label{observable1} \\
&&\;\;\int_0^t\!\!\int_\rr \gamma(s,x)\,H_n(t)(\mathrm{d}s,\mathrm{d}x)=\sqrt{n} \int_0^t\!\!\int_\rr \bigg(\!\int_0^x\! \gamma(s,y)\,\mathrm{d}y\!\bigg) (\rho^{(n)}(s)(\mathrm{d}x)\!-\!\rho(s)(\mathrm{d}x))\,\mathrm{d}s \label{observable2}
\end{eqnarray}
for functions $\gamma\in C_0(\rr)\cap L^1(\rr)$ and 
\begin{equation}
\gamma\in C_0([0,t]\times\rr):\;\gamma(s,\cdot)\in L^1(\rr)\;\text{for Lebesgue a.e.}\;s\in[0,t], 
\end{equation}
respectively. 
 
\medskip

Our main result can be stated as follows. 

\begin{theorem}\label{thm:CLT}
Suppose that Assumption \ref{main_asmp} holds and consider the mild solution $G$ of the SPDE
\begin{equation}\label{SPDE}
G_t=\big(b(R)\,G\big)_x+\bigg(\frac{\sigma(R)^2}{2}\,G\bigg)_{xx}
+\sigma(R)\,R_x^{1/2}\,\dot{W},\quad G(0,\cdot)=\beta(F_\lambda(\cdot)),
\end{equation}
where $R$ is the unique generalized solution to the Cauchy problem \eqref{PME}, $\dot{W}$ is a space-time white noise and $\beta$ is a standard Brownian bridge independent of $\dot{W}$. More specifically, let $G$ be the random field defined by
\begin{equation}\label{mild_def}
\begin{split}
G(t,x)=\int_\rr \beta(F_\lambda(y))\,p(0,y;t,x)\,\mathrm{d}y
+\int_0^t \!\int_\rr \sigma(R(s,y))\,R_x(s,y)^{1/2}\,p(s,y;t,x)\,\mathrm{d}W(s,y), \\
(t,x)\in[0,\infty)\times\rr,
\end{split}
\end{equation}
where $p$ denotes the transition density of the solution to the martingale problem associated with the operators $b(R(t,\cdot))\,\frac{\mathrm{d}}{\mathrm{d}x}+\frac{\sigma(R(t,\cdot))^2}{2}\,\frac{\mathrm{d}^2}{\mathrm{d}x^2}$, $t\ge0$ and the double integral should be understood in the It\^o sense.

Then, one has the following convergences: 
\begin{enumerate}[(a)]
\item The $\Mf$-valued processes $G_n$, $n\in\nn$ tend in the finite-dimensional distribution sense to $t\mapsto G(t,x)\,\mathrm{d}x$. 
\item The processes $H_n$, $n\in\nn$ taking values in $\Mft$, $t>0$
converge in the finite-dimensional distribution sense to $t\mapsto G(s,x)\,\mathbf{1}_{[0,t]\times\rr}(s,x)\,\mathrm{d}s\,\mathrm{d}x$, also jointly with the processes in (a).
\end{enumerate}
\end{theorem}

\begin{rmk}
The result of Theorem \ref{thm:CLT} shows that the evolution of the capital distribution in a large equity market, in which the logarithmic capitalizations follow \eqref{ips_sde}, can be approximated by $t\mapsto R(t,\cdot)+n^{-1/2} G(t,\cdot)$ up to an error of $o(n^{-1/2})$. This suggests that, if one combines the unnormalized versions of the capital distribution curves as in \cite[Figure 5.1]{Fe} to a \textit{surface} whose height encodes the relative rank associated with any given logarithmic capitalization at any given time, that surface should resemble a typical realization of the random surface $R+n^{-1/2}G$ indexed by $[0,\infty)\times\rr$. Consequently, properties of large equity markets, as captured by the observables of \eqref{observable1}, \eqref{observable2}, can be accessed through the corresponding observables of the random surface $R+n^{-1/2}G$.    
\end{rmk}

\begin{rmk}\label{rmk:mild}
It is immediate from \eqref{mild_def} that the mild solution $G$ of the SPDE \eqref{SPDE} is a mean zero Gaussian process with a covariance of
\begin{equation}
\begin{split}
& \int_\rr \int_\rr \big(F_\lambda(\min(y_1,y_2))-F_\lambda(y_1)\,F_\lambda(y_2)\big)\,p(0,y_1;t_1,x_1)\,p(0,y_2;t_2,x_2)\,\mathrm{d}y_1\,\mathrm{d}y_2 \\
& +\int_0^{\min(t_1,t_2)} \int_\rr \sigma(R(s,y))^2\,R_x(s,y)\,p(s,y;t_1,x_1)\,p(s,y;t_2,x_2)\,\mathrm{d}y\,\mathrm{d}s
\end{split}
\end{equation}
between any $G(t_1,x_1)$, $G(t_2,x_2)$. 
\end{rmk}

\begin{rmk}
For constant $b$ and $\sigma$ (when the particles are independent) and a fixed $t\ge0$ the convergence of $G_n(t)$, $n\in\nn$ falls into the framework of \cite[Theorem 2.1]{dGM} (see also \cite[Corollary 3.9]{BL}). The topology used there is the weak topology on $L^1(\rr)$ and the result is established using the Central Limit Theorem in cotype 2 spaces. Due to the dependence between the particles in the general case we cannot use the same machinery and instead need to start by establishing the tightness of $G_n(t)$, $n\in\nn$ directly. For this reason, we chose to work with the space $\Mf$ rather than $L^1(\rr)$, as it admits a more amenable compactness criterion. 
\end{rmk}

\smallskip

In the course of the proof of Theorem \ref{thm:CLT} we obtain the first quantitative propagation of chaos result for the particle system \eqref{ips_sde}. The general propagation of chaos paradigm (see \cite{Sz2}) suggests that for large $n$ the weak solution of \eqref{ips_sde} should be close to the strong solution of 
\begin{equation}\label{bar_part}
\begin{split}
\mathrm{d}\bar{X}^{(n)}_i(t)=b\big(R(t,\bar{X}^{(n)}_i(t))\big)\,\mathrm{d}t+\sigma\big(R(t,\bar{X}^{(n)}_i(t))\big)\,\mathrm{d}B^{(n)}_i(t),\;\;\;\bar{X}^{(n)}_i(0)=X^{(n)}_i(0), \\
i=1,\,2,\,\ldots,\,n,
\end{split}
\end{equation} 
where $B^{(n)}_1,\,B^{(n)}_2,\,\ldots,\,B^{(n)}_n$ are the standard Brownian motions from \eqref{ips_sde}. We refer to the discussion following Proposition \ref{PME_grad} below for the existence of a unique strong solution of \eqref{bar_part}. Writing $\bar{\rho}^{(n)}(t):=\frac{1}{n}\sum_{i=1}^n \delta_{\bar{X}^{(n)}_i(t)}$, $t\ge0$ for the path of empirical measures associated with the i.i.d. particles $\bar{X}^{(n)}_1,\,\bar{X}^{(n)}_2,\,\ldots,\,\bar{X}^{(n)}_n$ we aim to compare $\rho^{(n)}(\cdot)$ to $\bar{\rho}^{(n)}(\cdot)$. As a notion of distance we introduce for $p\ge1$ the Wasserstein metric $W_p$ on the space of probability measures on $\rr$ with finite moments up to order $p$:  
\begin{equation}\label{def:Wp}
W_p(\mu,\nu)=\inf_{(Y_1,Y_2)} \E\big[|Y_1-Y_2|^p\big]^{1/p},
\end{equation}
where the infimum is taken over all random vectors $(Y_1,Y_2)$ such that $Y_1$ is distributed according to $\mu$ and $Y_2$ according to $\nu$. Our quantitative propagation of chaos result then reads as follows.

\begin{theorem}\label{main1}
Suppose that Assumption \ref{main_asmp} holds. Then, for all $p>0$ and $T>0$ there exists a constant $C=C(p,T)<\infty$ such that
\begin{equation}\label{eq_poc}
\forall\,n\in\nn,\;1\le i\le n:\quad \E\Big[\sup_{0\le t\le T} \big|X^{(n)}_i(t)-\bar{X}^{(n)}_i(t)\big|^p\Big]\le C\,n^{-p/2}.
\end{equation}
In particular, when $p\ge1$ one has 
\begin{equation}\label{eq_poc'}
\forall\,n\in\nn:\quad\E\Big[\sup_{0\le t\le T} W_p\big(\rho^{(n)}(t),\bar{\rho}^{(n)}(t)\big)^p\Big]\le C\,n^{-p/2}.
\end{equation}  
\end{theorem}

\smallskip

The rest of the paper is structured as follows. In Section \ref{prelim} we prepare various results that are used in the proofs of Theorems \ref{thm:CLT} and \ref{main1}: some properties of Wasserstein distances and relations of the latter to empirical measures (from \cite{BL} and \cite{dGM}), as well as a PDE estimate for the solution of \eqref{PME} (from \cite{Gi}) and its implications for the associated diffusion process (including Gaussian lower and upper bounds on the transition density based on the results in \cite{Ar} and \cite{Kr2}). In Section \ref{sec:poc} we prove Theorem \ref{main1} by reducing it to the estimate of \cite[Theorem 4.8]{BL} on the expected Wasserstein distance between the empirical measure of an i.i.d. sample from the uniform distribution and the uniform distribution itself. Theorem \ref{main1} is then used in Section \ref{sec:tight} to establish the tightness of the finite-dimensional distributions of the processes $G_n$, $n\in\nn$ and $H_n$, $n\in\nn$ via a representation of $W_1$ for probability measures on $\rr$ in terms of their cumulative distribution functions. In Section \ref{sec:limit} we conclude the proof of Theorem \ref{thm:CLT} by
identifying the limit points of the finite-dimensional distributions of $G_n$, $n\in\nn$ and $H_n$, $n\in\nn$. Our argument relies on a prelimit version of the martingale problem associated with the SPDE \eqref{SPDE} (see Lemma \ref{lemma:testfunction}) and an appropriate coupling construction (see the proof of Proposition \ref{prop:full_iden}).  

\medskip

\noindent\textbf{Acknowledgements.} We would like to thank Cameron Bruggeman for an enlightening discussion at an early stage of the preparation of this paper. We also thank Ioannis Karatzas for his many helpful comments. 

\section{Preliminaries} \label{prelim}

\subsection{Wasserstein distances and empirical measures}

For $p\ge1$ consider two probability measures $\mu$, $\nu$ on $\rr$ having finite moments up to order $p$. Let $F_\mu$, $F_\nu$ be their cumulative distribution functions and $q_\mu$, $q_\nu$ be their quantile functions. The following well-known representations of $W_p(\mu,\nu)$ (see e.g. \cite[Section 2.3]{BL}) are used repeatedly below. 

\begin{proposition}
In the setup of the preceding paragraph it holds
\begin{eqnarray}
&& W_1(\mu,\nu)=\int_\rr \big|F_\mu(x)-F_\nu(x)\big|\,\mathrm{d}x, \label{W1repr} \\
&& W_p(\mu,\nu)=\bigg(\int_0^1 \big|q_\mu(a)-q_\nu(a)\big|^p\,\mathrm{d}a\bigg)^{1/p},\quad p\ge1. \label{quantile_repr}
\end{eqnarray} 
\end{proposition}

\smallskip

In addition, we prepare estimates on the expected Wasserstein distances between the empirical measure of an i.i.d. sample from the uniform distribution and the uniform distribution itself. These are taken from \cite[Theorem 4.8]{BL}.
\begin{proposition}\label{prop:unif_sample}
Let $U_1,\,U_2,\,\ldots$ be i.i.d. according to the uniform distribution $\upsilon$ on $[0,1]$. Then, there exists a constant $C<\infty$ such that 
\begin{equation}\label{eq:unif_sample}
\E\bigg[W_p\bigg(\frac{1}{n}\sum_{i=1}^n \delta_{U_i},\upsilon\bigg)^p\bigg]^{1/p}\le C\,p^{1/2}\,n^{-1/2},\quad p\ge1,\quad n\in\nn.
\end{equation}
\end{proposition}

\smallskip

Finally, we recall the Functional Central Limit Theorem for empirical cumulative distribution functions from \cite[Theorem 2.1]{dGM} (see also \cite[Corollary 3.9 and the discussion of the functional $J_1$ on p. 25]{BL}). This result gives rise to the initial condition in \eqref{SPDE}.

\begin{proposition}\label{CLT_dGM}
Let Assumption \ref{main_asmp}(a) be satisfied. Then, the sequence $G_n(0,\cdot)$, $n\in\nn$ converges in law weakly in $L^1(\rr)$ (and therefore in $\Mf$) to $\beta(F_\lambda(\cdot))$, where $\beta$ is a standard Brownian bridge. 
\end{proposition}

\subsection{Porous medium equation and associated diffusion process}\label{sec:PMEetc}

We turn to the properties of the generalized solution to the problem \eqref{PME} and the associated diffusion process. First, we briefly recall \cite[Definition 3]{Gi} of such a generalized solution (see also the original reference \cite[Definition 1.1]{DK}).

\begin{definition}\label{Gi_def}
A bounded continuous nonnegative function $R$ with $R(0,\cdot)=F_\lambda(\cdot)$ is called a generalized solution of the Cauchy problem \eqref{PME} if 
\begin{equation}
\begin{split}
\int_{t_1}^{t_2} \!\!\int_{x_1}^{x_2} \!\! \zeta_x\,B(R)+\zeta_{xx}\,\Sigma(R)+\zeta_t\,R\,\mathrm{d}x\,\mathrm{d}t \!=\!\! \int_{x_1}^{x_2} \!\! \zeta(t_2,\cdot)\,R(t_2,\cdot)\,\mathrm{d}x \!-\! \int_{x_1}^{x_2} \!\! \zeta(t_1,\cdot)\,R(t_1,\cdot)\,\mathrm{d}x \\
+ \int_{t_1}^{t_2} \zeta_x(t,x_2)\,\Sigma(R(t,x_2))\,\mathrm{d}t - \int_{t_1}^{t_2} \zeta_x(t,x_1)\,\Sigma(R(t,x_1))\,\mathrm{d}t
\end{split}
\end{equation}
for all $0\le t_1<t_2$, $x_1<x_2$ and functions $\zeta:\,[t_1,t_2]\times[x_1,x_2]\to\rr$ which are continuously differentiable in $t$, twice continuously differentiable in $x$ and satisfy $\zeta(\cdot,x_1)=\zeta(\cdot,x_2)=0$. 
\end{definition}

In view of Assumption \ref{main_asmp} and $\min_{a\in[0,1]} \Sigma'(a)=\min_{a\in[0,1]} \frac{1}{2}\,\sigma(a)^2>0$, we can combine \cite[Theorems 4 and 7]{Gi} to the following proposition. 

\begin{proposition}\label{PME_grad}
Let Assumption \ref{main_asmp} be satisfied. Then, the Cauchy problem \eqref{PME} admits a unique generalized solution $R$. Moreover, its distributional derivative $R_x$ can be represented by a bounded function on any strip of the form $[0,T]\times\rr$.  
\end{proposition}

\smallskip

We conclude the subsection with a discussion of the SDE  
\begin{equation}\label{PME_SDE}
\mathrm{d}\bar{X}(t)=b\big(R(t,\bar{X}(t))\big)\,\mathrm{d}t+\sigma\big(R(t,\bar{X}(t))\big)\,\mathrm{d}B(t)
\end{equation}
satisfied by each of the processes $\bar{X}^{(n)}_i$. Assumption \ref{main_asmp} and Proposition \ref{PME_grad} guarantee that the functions $x\mapsto b\big(R(t,x)\big)$ and $x\mapsto \sigma\big(R(t,x)\big)$ are Lipschitz with uniformly bounded Lipschitz constants on every compact interval of $t$'s. Consequently, there exists a unique strong solution of \eqref{PME_SDE} for the initial condition $\lambda$ of Assumption \ref{main_asmp} or any deterministic initial condition (see e.g. \cite[Chapter 5, Theorems 2.5 and 2.9]{KS}). In addition, $\bar{X}$ is the unique solution of the martingale problem associated with the operators $b(R(t,\cdot))\,\frac{\mathrm{d}}{\mathrm{d}x}+\frac{\sigma(R(t,\cdot))^2}{2}\,\frac{\mathrm{d}^2}{\mathrm{d}x^2}$, $t\ge0$ and therefore a strong Markov process (see \cite[Theorems 7.2.1 and 6.2.2]{SV}). For the initial condition $\lambda$, Assumption \ref{main_asmp} allows us to apply \cite[Corollary 1.13]{JR} to identify the one-dimensional distributions of the solution to the nonlinear martingale problem therein with $\rho(t)$, $t\ge0$, so that the solution itself is given by the law ${\mathcal L}(\bar{X})$ of $\bar{X}$ and therefore
\begin{equation}\label{lawbarX}
{\mathcal L}(\bar{X}(t))=\rho(t),\quad t\ge0.
\end{equation}

\smallskip

We now aim to apply the results of \cite{Ar} to conclude that under Assumption \ref{main_asmp} the transition density of $\bar{X}$ exists and satisfies Gaussian lower and upper bounds. To identify the transition density of $\bar{X}$ with the weak fundamental solution of a parabolic PDE as in \cite[Theorem 5]{Ar} we fix a $T>0$ and consider the Cauchy problem
\begin{equation}\label{Krylov_PDE}
u_t+b(R)\,u_x+\frac{\sigma(R)^2}{2}\,u_{xx}=f,\quad u(T,\cdot)=0,
\end{equation}
where $f\in L^2([0,T]\times\rr)\cap L^\infty([0,T]\times\rr)$. We note that $b(R)$ and $\frac{\sigma(R)^2}{2}$ are bounded and that $x\mapsto \frac{\sigma(R(t,x))^2}{2}$ are Lipschitz with uniformly bounded Lipschitz constants for $t\in[0,T]$. Hence, according to \cite[Theorem 2.1 and Remark 2.2]{Kr2} there exists a unique solution $u$ of \eqref{Krylov_PDE} with $u,u_t,u_x,u_{xx}\in L^2([0,T]\times\rr)$ and it is given by
\begin{equation}\label{Feynman-Kac}
u(t,x)=-\E\bigg[\int_t^T f(r,\bar{X}(r))\,\mathrm{d}r\,\bigg|\,\bar{X}(t)=x\bigg],\quad (t,x)\in[0,T]\times\rr.
\end{equation}
In particular, $u\in L^\infty([0,T]\times\rr)$, so that with $g(t,x):=-f(T-t,x)$, $S(t,x):=R(T-t,x)$, $(t,x)\in[0,T]\times\rr$ the function $v(t,x):=u(T-t,x)$, $(t,x)\in[0,T]\times\rr$ is a weak solution of 
\begin{equation}\label{backward_div}
v_t-\big(b(S)-\sigma(S)\,\sigma'(S)\,S_x\big)\,v_x-\bigg(\frac{\sigma(S)^2}{2}\,v_x\bigg)_x=g,\quad v(0,\cdot)=0
\end{equation}
in the sense of \cite[Theorem 5(ii)]{Ar}. The latter theorem is applicable, since $b(S)-\sigma(S)\,\sigma'(S)\,S_x$ and $\frac{\sigma(S)^2}{2}$ are bounded on $[0,T]\times\rr$ and $\frac{\sigma(S)^2}{2}$ is bounded away from $0$ on $[0,T]\times\rr$ by Assumption \ref{main_asmp} and Proposition \ref{PME_grad}. Comparing the conclusion of \cite[Theorem 5(ii)]{Ar} with \eqref{Feynman-Kac} we obtain the existence of the transition density $p(t,x;r,z)$ of $\bar{X}$ and recognize $p(T-r,x;T-t,z)$ as the weak fundamental solution corresponding to the PDE in \eqref{backward_div}. Thus, \cite[Theorem 10(ii)]{Ar} yields the following result. 

\begin{proposition}\label{prop:heat_kernel}
Let Assumption \ref{main_asmp} be satisfied. Then, the process $\bar{X}$ has a transition density $p$ such that
\begin{equation}\label{eq:heat_kernel}
\begin{split}
\forall\,T>0:\;\; C^{-1} (r-t)^{-1/2}\,e^{-C(z-x)^2/(r-t)}\le p(t,x;r,z)\le C (r-t)^{-1/2}\,e^{-C^{-1}(z-x)^2/(r-t)}, \\
0\le t<r\le T,\;x,z\in\rr
\end{split}
\end{equation}
with $C\in(1,\infty)$ possibly depending on $T$. In particular, if $\bar{X}(0)$ is distributed according to $\lambda$, then
\begin{equation}\label{moment_est}
\forall\,T>0:\quad\sup_{0\le t\le T} \E\big[|\bar{X}(t)|^{2+\eta}\big]<\infty.
\end{equation}
\end{proposition}

\section{Propagation of chaos estimates} \label{sec:poc}

This section is devoted to the proof of Theorem \ref{main1}. 

\medskip

\noindent\textbf{Proof of Theorem \ref{main1}. Step 1.} Fix any $p\ge2$ and $T>0$. We aim to employ Proposition \ref{prop:unif_sample} and to do so we are going to estimate the left-hand side of \eqref{eq_poc} by a quantity involving the left-hand side of \eqref{eq:unif_sample}. To this end, we first observe that the pairs $(X^{(n)}_i,\bar{X}^{(n)}_i)$, $i=1,\,2,\,\ldots,\,n$ have the same distribution (due to the weak uniqueness for \eqref{ips_sde} and the strong uniqueness for \eqref{bar_part}) and therefore the left-hand side of \eqref{eq_poc} can be rewritten in the symmetrized form
\begin{equation}\label{poc_symmetrized}
\frac{1}{n}\,\sum_{i=1}^n\,\E\Big[\sup_{0\le t\le T} \big|X^{(n)}_i(t)-\bar{X}^{(n)}_i(t)\big|^p\Big].
\end{equation}
 
Next, we use the SDEs \eqref{ips_sde} and \eqref{bar_part} satisfied by $X^{(n)}_i$ and $\bar{X}^{(n)}_i$, the elementary inequality  
\begin{equation}\label{ineq_elem_p}
(r_1+r_2)^p\le 2^{p-1}(r_1^p+r_2^p),\quad r_1,r_2\ge 0, 
\end{equation}
the Burkholder-Davis-Gundy inequality (see e.g. \cite[Chapter 3, Theorem 3.28]{KS}) and the Lipschitz property of $b$ and $\sigma$ to find for all $t\in[0,T]$ and $i=1,\,2,\,\ldots,\,n$:
\begin{equation}\label{eq_poc_0}
\begin{split}
\E\Big[\sup_{0\le s\le t} \big|X^{(n)}_i(s)-\bar{X}^{(n)}_i(s)|^p\Big]\le &\,C\,\E\bigg[\bigg(\!\int_0^t\!\big|F_{\rho^{(n)}(s)}(X^{(n)}_i(s))\!-\!R(s,\bar{X}^{(n)}_i(s))\big|\,\mathrm{d}s\!\bigg)^p\bigg] \\
&+\!C\,\E\bigg[\bigg(\!\int_0^t\!\big|F_{\rho^{(n)}(s)}(X^{(n)}_i(s))\!-\!R(s,\bar{X}^{(n)}_i(s))\big|^2\,\mathrm{d}s\!\bigg)^{p/2}\bigg],
\end{split}
\end{equation}
where $C<\infty$ depends only on $p$ and the Lipschitz constants of $b$ and $\sigma$. Applying Jensen's inequality to each of the summands on the right-hand side of \eqref{eq_poc_0} we obtain the further upper bound
\begin{equation}
C\,\E\bigg[\int_0^t \big|F_{\rho^{(n)}(s)}(X^{(n)}_i(s))-R(s,\bar{X}^{(n)}_i(s))\big|^p\,\mathrm{d}s\bigg],
\end{equation}
where $C<\infty$ can be chosen in terms of $T$, $p$ and the Lipschitz constants of $b$ and $\sigma$. 

\medskip

Another application of \eqref{ineq_elem_p} gives
\begin{equation}\label{eq_poc_1}
\begin{split}
& C\,\E\bigg[\int_0^t \big|F_{\rho^{(n)}(s)}(X^{(n)}_i(s))-R(s,X^{(n)}_i(s))\big|^p\,\mathrm{d}s\bigg] \\
& +C\,\E\bigg[\int_0^t \big|R(s,X^{(n)}_i(s))-R(s,\bar{X}^{(n)}_i(s))\big|^p\,\mathrm{d}s\bigg],
\end{split}
\end{equation}
where $C<\infty$ is still a function of $T$, $p$ and the Lipschitz constants of $b$ and $\sigma$ only. Now, we take the average of the first summands in \eqref{eq_poc_1} over $i=1,\,2,\,\ldots,\,n$ and get
\begin{equation}\label{eq_poc_2}
\begin{split}
& \frac{C}{n}\,\sum_{i=1}^n\,\E\bigg[\int_0^t \big|F_{\rho^{(n)}(s)}(X^{(n)}_i(s))-R(s,X^{(n)}_i(s))\big|^p\,\mathrm{d}s\bigg] \\
& =C \int_0^t \E\bigg[\frac{1}{n}\,\sum_{i=1}^n \big|F_{\rho^{(n)}(s)}(X^{(n)}_i(s))-R(s,X^{(n)}_i(s))\big|^p\bigg]\,\mathrm{d}s \\
& =C \int_0^t \E\bigg[\frac{1}{n}\,\sum_{k=1}^n \big|F_{\rho^{(n)}(s)}(X^{(n)}_{(k)}(s))-R(s,X^{(n)}_{(k)}(s))\big|^p\bigg]\,\mathrm{d}s,
\end{split}
\end{equation}
where $X^{(n)}_{(1)}(s)\le X^{(n)}_{(2)}(s)\le\cdots\le X^{(n)}_{(n)}(s)$ are the order statistics of the vector $\big(X^{(n)}_1(s),X^{(n)}_2(s),\ldots,X^{(n)}_n(s)\big)$. 

\medskip

At this point, \cite[Theorem on p. 439]{Kr1} for the function $y\mapsto\sum_{1\le i<j\le n} \mathbf{1}_{\{y_i=y_j\}}$ on $\rr^n$ reveals that with probability one it holds $F_{\rho^{(n)}(s)}(X^{(n)}_{(k)}(s))=\frac{k}{n}$, $k=1,\,2,\,\ldots,\,n$ for Lebesgue a.e. $s\in[0,T]$. This and \eqref{ineq_elem_p} allow to estimate the end result of \eqref{eq_poc_2} from above by 
\begin{equation}\label{eq_poc_3}
\begin{split}
C \int_0^t \E\bigg[\frac{1}{n}\,\sum_{k=1}^n \bigg|\frac{k}{n}-R(s,\bar{X}^{(n)}_{(k)}(s))\bigg|^p\bigg]
+\E\bigg[\frac{1}{n}\,\sum_{k=1}^n \big|R(s,\bar{X}^{(n)}_{(k)}(s))-R(s,X^{(n)}_{(k)}(s))\big|^p\bigg]\,\mathrm{d}s, 
\end{split}
\end{equation}
where $\bar{X}^{(n)}_{(1)}(s)\le \bar{X}^{(n)}_{(2)}(s)\le\cdots\le \bar{X}^{(n)}_{(n)}(s)$ are the order statistics of the vector $\big(\bar{X}^{(n)}_1(s),\bar{X}^{(n)}_2(s),\ldots,\bar{X}^{(n)}_n(s)\big)$ and $C<\infty$ depends on $T$, $p$ and the Lipschitz constants of $b$ and $\sigma$ only. 

\medskip

\noindent\textbf{Step 2.} Relying on the representation \eqref{quantile_repr} we readily identify the quantity $\E\big[\frac{1}{n}\sum_{k=1}^n \big|\frac{k}{n}-R(s,\bar{X}^{(n)}_{(k)}(s))\big|^p\big]$ in \eqref{eq_poc_3} as 
\begin{equation}\label{eq_poc_4}
\E\bigg[W_p\bigg(\frac{1}{n}\sum_{k=1}^n \delta_{k/n},\,\frac{1}{n}\sum_{k=1}^n \delta_{R(s,\bar{X}^{(n)}_{(k)}(s))}\bigg)^p\bigg].
\end{equation}
The observation \eqref{lawbarX} reveals $R(s,\bar{X}^{(n)}_{(1)}(s))\le R(s,\bar{X}^{(n)}_{(2)}(s))\le\cdots\le R(s,\bar{X}^{(n)}_{(n)}(s))$ as the order statistics of an i.i.d. sample from the uniform distribution on $[0,1]$. This, the triangle inequality for $W_p$ and  \eqref{ineq_elem_p} imply that the expectation in \eqref{eq_poc_4} is bounded above by
\begin{equation}\label{eq_poc_5}
2^{p-1}\,\E\bigg[W_p\bigg(\frac{1}{n}\sum_{k=1}^n \delta_{k/n},\,\upsilon\bigg)^p\bigg]+2^{p-1}\,\E\bigg[W_p\bigg(\upsilon,\,\frac{1}{n}\sum_{i=1}^n \delta_{U_i}\bigg)^p\bigg]
\end{equation}
in the notation of Proposition \ref{prop:unif_sample}. Using the representation \eqref{quantile_repr} for the first expectation in \eqref{eq_poc_5} and Proposition \ref{prop:unif_sample} for the second expectation in \eqref{eq_poc_5} we end up with the upper bound 
\begin{equation}\label{eq_poc_6}
2^{p-1}\,n^{-p}+2^{p-1}\,C^p\,p^{p/2}\,n^{-p/2},
\end{equation} 
where $C$ is the constant in Proposition \ref{prop:unif_sample}.

\medskip

\noindent\textbf{Step 3.} Putting the estimates \eqref{eq_poc_1}, \eqref{eq_poc_3} and \eqref{eq_poc_6} together we arrive at the inequality 
\begin{equation}\label{eq_poc_7}
\begin{split}
&\frac{1}{n}\,\sum_{i=1}^n\,\E\Big[\sup_{0\le s\le t} \big|X^{(n)}_i(s)-\bar{X}^{(n)}_i(s)\big|^p\Big] \\
&\le C\,\int_0^t \bigg(n^{-p}+n^{-p/2}+\E\bigg[\frac{1}{n}\,\sum_{k=1}^n\big|R(s,X^{(n)}_{(k)}(s))-R(s,\bar{X}^{(n)}_{(k)}(s))\big|^p\bigg] \\
&\qquad\qquad\qquad\qquad\quad\;\;\; +\E\bigg[\frac{1}{n}\,\sum_{i=1}^n\big|R(s,X^{(n)}_i(s))-R(s,\bar{X}^{(n)}_i(s))\big|^p\bigg]\bigg)\,\mathrm{d}s
\end{split}
\end{equation}
for all $t\in[0,T]$, where $C<\infty$ is a function of $T$, $p$ and the Lipschitz constants of $b$ and $\sigma$. Moreover, the functions $x\mapsto R(s,x)$ are Lipschitz with uniformly bounded Lipschitz constants as $s$ varies in $[0,T]$ by Proposition \ref{PME_grad} and 
\begin{equation}
\frac{1}{n}\,\sum_{k=1}^n \big|X^{(n)}_{(k)}(s)-\bar{X}^{(n)}_{(k)}(s)\big|^p=W_p\big(\rho^{(n)}(s),\bar{\rho}^{(n)}(s)\big)^p\le \frac{1}{n}\,\sum_{i=1}^n \big|X^{(n)}_i(s)-\bar{X}^{(n)}_i(s)\big|^p
\end{equation}
by the representation \eqref{quantile_repr} and the definition of $W_p$ in \eqref{def:Wp}, so that for all $t\in[0,T]$:
\begin{equation}\label{313}
\begin{split}
\frac{1}{n}\sum_{i=1}^n \E\Big[\sup_{0\le s\le t} \big|X^{(n)}_i(s)-\bar{X}^{(n)}_i(s)\big|^p\Big]\le &\,C\big(n^{-p}+n^{-p/2}\big)t \\
&+\!C\!\int_0^t \frac{1}{n}\sum_{i=1}^n \E\Big[\sup_{0\le r\le s} \big|X^{(n)}_i(r)-\bar{X}^{(n)}_i(r)\big|^p\Big]\,\mathrm{d}s,
\end{split}
\end{equation}
where $C<\infty$ depends on $T$, $p$, the Lipschitz constants of $b$ and $\sigma$ and the supremum of $R_x$ on $[0,T]\times\rr$ only. The desired estimate \eqref{eq_poc} is a consequence of \eqref{313} due to the representation \eqref{poc_symmetrized} and Gronwall's lemma. 

\medskip

\noindent\textbf{Step 4.} For $p\in(0,2)$, we choose a $p'\in[2,\infty)$ and deduce \eqref{eq_poc} for $p$ from \eqref{eq_poc} for $p'$ by means of the inequality
\begin{equation}
\frac{1}{n}\sum_{i=1}^n \E\Big[\sup_{0\le t\le T} \big|X^{(n)}_i(t)-\bar{X}^{(n)}_i(t)\big|^p\Big]\le \bigg(\frac{1}{n}\sum_{i=1}^n \E\Big[\sup_{0\le t\le T} \big|X^{(n)}_i(t)-\bar{X}^{(n)}_i(t)\big|^{p'}\Big]\bigg)^{p/p'}.
\end{equation} 
Finally, we obtain \eqref{eq_poc'} from \eqref{eq_poc} via the chain of estimates
\begin{equation}
\begin{split}
& \E\Big[\sup_{0\le t\le T} W_p\big(\rho^{(n)}(t),\bar{\rho}^{(n)}(t)\big)\Big]^p \le \E\Big[\sup_{0\le t\le T} W_p\big(\rho^{(n)}(t),\bar{\rho}^{(n)}(t)\big)^p\Big] \\
& \le\E\bigg[\sup_{0\le t\le T}\,\frac{1}{n}\,\sum_{i=1}^n\, \big|X^{(n)}_i(t)-\bar{X}^{(n)}_i(t)\big|^p\bigg]\le 
\frac{1}{n}\,\sum_{i=1}^n\,\E\Big[\sup_{0\le t\le T} \big|X^{(n)}_i(t)-\bar{X}^{(n)}_i(t)\big|^p\Big]
\end{split}
\end{equation}
valid for all $p\ge1$. \ep

\section{Existence of subsequential limits} \label{sec:tight}

The main result of this section is the next proposition establishing the existence of subsequential limits for the finite-dimensional distributions of the fluctuation processes $G_n$, $n\in\nn$ and $H_n$, $n\in\nn$. It serves as a key ingredient in the proof of Theorem \ref{thm:CLT}.

\begin{proposition}\label{prop:tight}
Suppose that Assumption \ref{main_asmp} is satisfied. Then, for all $m\in\nn$ and $0<t_1<\cdots<t_m$ every subsequence of 
\begin{equation}\label{FDDseq}
\big(G_n(0), G_n(t_1), \ldots, G_n(t_m), H_n(t_1), H_n(t_2), \ldots, H_n(t_m)\big),\quad n\in\nn
\end{equation}
has a further subsequence which converges in law in
\begin{equation*}
\Mf^{m+1}\times M_{\mathrm{fin}}([0,t_1]\times\rr)\times
M_{\mathrm{fin}}([0,t_2]\times\rr)\times\cdots\times M_{\mathrm{fin}}([0,t_m]\times\rr).
\end{equation*}
\end{proposition}

\smallskip

\noindent\textbf{Proof.} By Prokhorov's Theorem in the form of \cite[Corollary on p. 119]{FGH} it suffices to show that the laws of the random vectors in \eqref{FDDseq} form a uniformly tight sequence. Moreover, since products of compact sets are compact, we only need to prove that for all $s\ge0$ and $t>0$ the laws associated with the sequences $G_n(s)$, $n\in\nn$ and $H_n(t)$, $n\in\nn$ are uniformly tight. In view of the Banach-Alaoglu Theorem (see e.g. \cite[Chapter 12, Theorem 3]{La}), this is the case for any fixed $s\ge0$ and $t>0$ if for all $\epsilon>0$ there exists a $C_\epsilon<\infty$ such that
\begin{equation}\label{TVform}
\forall\,n\in\nn:\quad\pp\big(\|G_n(s)\|_{TV}>C_\epsilon\big)<\epsilon\quad\text{and}\quad \pp\big(\|H_n(t)\|_{TV}>C_\epsilon\big)<\epsilon,
\end{equation}
where $\|\cdot\|_{TV}$ stands for the total variation norm.

\medskip

By the definitions of $G_n(s)$, $n\in\nn$ and $H_n(t)$, $n\in\nn$ in \eqref{def:Gn} and \eqref{def:Hn} the two inequalities of \eqref{TVform} can be rewritten as
\begin{eqnarray}
&& \pp\bigg(\sqrt{n}\,\int_\rr \big|F_{\rho^{(n)}(s)}(x)-R(s,x)\big|\,\mathrm{d}x>C_\epsilon\bigg)<\epsilon, \\ 
&& \pp\bigg(\sqrt{n}\,\int_0^t \int_\rr \big|F_{\rho^{(n)}(r)}(x)-R(r,x)\big|\,\mathrm{d}x\,\mathrm{d}r>C_\epsilon\bigg)<\epsilon.
\end{eqnarray}
The representation \eqref{W1repr} allows to rewrite these further as 
\begin{equation}\label{W1form}
\pp\Big(\sqrt{n}\,W_1(\rho^{(n)}(s),\rho(s))>C_\epsilon\Big)
<\epsilon,\quad
\pp\bigg(\sqrt{n}\,\int_0^t W_1(\rho^{(n)}(r),\rho(r))\,\mathrm{d}r>C_\epsilon\bigg)<\epsilon.
\end{equation}
Applying Markov's inequality, the triangle inequality for $W_1$ and Fubini's Theorem we bound the two probabilities in \eqref{W1form} from above by
\begin{eqnarray}
&& \frac{\sqrt{n}}{C_\epsilon}\,\E\big[W_1(\rho^{(n)}(s),\bar{\rho}^{(n)}(s))\big]+\frac{\sqrt{n}}{C_\epsilon}\,\E\big[W_1(\bar{\rho}^{(n)}(s),\rho(s))\big], \label{W1fixedtime} \\
&& \frac{\sqrt{n}}{C_\epsilon}\, 
\E\bigg[\int_0^t W_1(\rho^{(n)}(r),\bar{\rho}^{(n)}(r))\,\mathrm{d}r\bigg] + \frac{\sqrt{n}}{C_\epsilon}\,\int_0^t 
\E\big[W_1(\bar{\rho}^{(n)}(r),\rho(r))\big]\,\mathrm{d}r, \label{W1integrated}
\end{eqnarray} 
respectively. In view of \eqref{eq_poc'}, \cite[Theorem 3.2 and the discussion of the functional $J_1$ on p. 25]{BL}, \eqref{moment_est} and \eqref{lawbarX}, we can make the estimates \eqref{W1fixedtime}, \eqref{W1integrated} smaller than $\epsilon$ for all $n\in\nn$ by choosing a large enough $C_\epsilon<\infty$. \ep 

\section{Identification of subsequential limits}\label{sec:limit} 

In this section we identify the subsequential limits of Proposition \ref{prop:tight} and complete the proof of Theorem \ref{thm:CLT}. The next proposition is the first step towards such an identification. 

\begin{proposition}\label{prop:Gaussians}
Suppose that Assumption \ref{main_asmp} holds and let 
\begin{equation}
\big(G_\infty(0), G_\infty(t_1), \ldots, G_\infty(t_m), H_\infty(t_1), H_\infty(t_2), \ldots, H_\infty(t_m)\big)
\end{equation}
be a limit point in law of the sequence in \eqref{FDDseq}. Then, the joint distribution of
\begin{equation}\label{limitGaussians}
\begin{split}
&\int_\rr \gamma(t_\ell,x)\,G_\infty(t_\ell)(\mathrm{d}x)
-\int_\rr \gamma(0,x)\,G_\infty(0)(\mathrm{d}x) \\
& -\int_0^{t_\ell} \int_\rr \Big(\gamma_s(s,x)+\gamma_x(s,x)\, b(R(s,x))+\gamma_{xx}(s,x)\,\frac{\sigma(R(s,x))^2}{2}\Big)\,H_\infty(t_\ell)(\mathrm{d}s,\mathrm{d}x), \\
&\int_\rr \gamma(0,x)\,G_\infty(0)(\mathrm{d}x),
\end{split}
\end{equation}
as $\ell$ and $\gamma$ vary over $\{1,2,\ldots,m\}$ and the space of functions on $[0,t_\ell]\times\rr$ which are continuously differentiable in $s$, twice continuously differentiable in $x$ and compactly supported, coincides with that of
\begin{equation}\label{limitGaussians'}
\int_0^{t_\ell} \int_\rr \gamma(s,x)\,\sigma(R(s,x))\,R_x(s,x)^{1/2} \,\mathrm{d}W(s,x),\quad \int_\rr \gamma(0,x)\,\beta(F_\lambda(x))\,\mathrm{d}x
\end{equation}
in the notation of Theorem \ref{thm:CLT}.
\end{proposition}

The proof of Proposition \ref{prop:Gaussians} relies on a suitable prelimit version of its statement. For every fixed $n\in\nn$ let $B_n$, $\Sigma_n$ be the piecewise constant functions on $[0,1]$ with jumps at $\frac{1}{n},\,\frac{2}{n},\,\ldots,\,1$ and  
\begin{equation}\label{BnSn}
B_n(k/n)=\frac{1}{n}\,\sum_{j=1}^k b(j/n),\;\;  
\Sigma_n(k/n)=\frac{1}{n}\,\sum_{j=1}^k \frac{\sigma(j/n)^2}{2},\quad k=0,\,1,\,\ldots,\,n.
\end{equation}

\begin{lemma}\label{lemma:testfunction}
Suppose that Assumption \ref{main_asmp} is satisfied. Then, for any $n\in\nn$, $t>0$ and function $\gamma$ on $[0,t]\times\rr$ which is continuously differentiable in $s$, twice continuously differentiable in $x$ and compactly supported it holds 
\begin{equation}\label{eq:main_repr}
\begin{split}
& \int_\rr \gamma(t,x)\,G_n(t)(\mathrm{d}x)
-\int_\rr \gamma(0,x)\,G_n(0)(\mathrm{d}x) \\
&-\int_0^t \int_\rr \int_0^1 \bigg(\gamma_s(s,x)+\gamma_x(s,x)\,b(a F_{\rho^{(n)}(s)}(x)+(1-a)R(s,x)) \\
&\qquad\qquad\qquad\,
+\gamma_{xx}(s,x)\,\frac{\sigma(a F_{\rho^{(n)}(s)}(x)+(1-a) R(s,x))^2}{2}\bigg)\,\mathrm{d}a\,\,H_n(t)(\mathrm{d}s,\mathrm{d}x) \\
&=-\frac{1}{\sqrt{n}}\,\sum_{i=1}^n \int_0^t \gamma(s,X^{(n)}_i(s))\,\sigma(F_{\rho^{(n)}(s)}(X^{(n)}_i(s)))\,\mathrm{d}B^{(n)}_i(s) \\
&\quad\,+\sqrt{n}\,\int_0^t \int_\rr \big(\gamma_x(s,x)\,(B_n-B)(F_{\rho^{(n)}(s)}(x))+\gamma_{xx}(s,x)\,(\Sigma_n-\Sigma)(F_{\rho^{(n)}(s)}(x))\big)\,\mathrm{d}x\,\mathrm{d}s.
\end{split}
\end{equation}
\end{lemma}

\smallskip

\noindent\textbf{Proof of Lemma \ref{lemma:testfunction}.} Fixing $n$, $t$ and $\gamma$ as described we observe that Definition \ref{Gi_def} of a generalized solution to the Cauchy problem \eqref{PME} implies 
\begin{equation}\label{PMEh}
\begin{split}
& \int_\rr \gamma(t,x)\,R(t,x)\,\mathrm{d}x-\int_\rr \gamma(0,x)\,R(0,x)\,\mathrm{d}x \\
&=\int_0^t \int_\rr \gamma_s(s,x)\,R(s,x)+\gamma_x(s,x)\,B(R(s,x))+\gamma_{xx}(s,x)\,\Sigma(R(s,x))\,\mathrm{d}x\,\mathrm{d}s.
\end{split}
\end{equation}
To find a version of the identity \eqref{PMEh} with $F_{\rho^{(n)}(\cdot)}(\cdot)$ in place of $R(\cdot,\cdot)$ we apply It\^o's formula for $\Gamma(s,x):=-\int_x^\infty \gamma(s,y)\,\mathrm{d}y$ and obtain 
\begin{equation}\label{eq_expansion1_26}
\begin{split}
&\int_\rr \Gamma(t,x)\,\rho^{(n)}(t)(\mathrm{d}x)-\int_\rr \Gamma(0,x)\,\rho^{(n)}(0)(\mathrm{d}x) \\
&=\frac{1}{n}\,\sum_{i=1}^n \int_0^t \gamma(s,X^{(n)}_i(s))\,\sigma(F_{\rho^{(n)}(s)}(X^{(n)}_i(s)))\,\mathrm{d}B^{(n)}_i(s) \\
&\quad+\int_0^t\!\int_\rr\bigg(\Gamma_s(s,x)\!+\!\Gamma_x(s,x) b(F_{\rho^{(n)}(s)}(x))
\!+\!\Gamma_{xx}(s,x) \frac{\sigma(F_{\rho^{(n)}(s)}(x))^2}{2}\bigg)\rho^{(n)}(s)(\mathrm{d}x)\,\mathrm{d}s.
\end{split}
\end{equation}

\smallskip

Next, we use summation by parts (note that $\lim_{x\to\infty} \Gamma(s,x)=0$ and $\lim_{x\to\infty} \Gamma_s(s,x)=0$ for all $s\in[0,t]$ by the compact support assumption on $\gamma$) to compute
\begin{eqnarray}
&&\quad\int_\rr \Gamma(s,x) \rho^{(n)}(s)(\mathrm{d}x)\!=\!-\!\int_\rr \gamma(s,x) F_{\rho^{(n)}(s)}(x)\,\mathrm{d}x,\,s\in\{0,t\}, \label{sbp1} \\
&&\quad\int_\rr \Gamma_s(s,x) \rho^{(n)}(s)(\mathrm{d}x)\!=\!-\!\int_\rr \gamma_s(s,x) F_{\rho^{(n)}(s)}(x)\,\mathrm{d}x,\,s\in[0,t], \label{sbp2} \\
&&\quad\int_\rr \Gamma_x(s,x) b(F_{\rho^{(n)}(s)}(x)) \rho^{(n)}(s)(\mathrm{d}x)\!=\!-\!\int_\rr \gamma_x(s,x) B_n(F_{\rho^{(n)}(s)}(x))\,\mathrm{d}x,\,s\in[0,t], \label{sbp3} \\
&&\quad\int_\rr \Gamma_{xx}(s,x) \frac{\sigma(F_{\rho^{(n)}(s)}(x))^2}{2} \rho^{(n)}(s)(\mathrm{d}x)\!=\!-\!\int_\rr \gamma_{xx}(s,x) \Sigma_n(F_{\rho^{(n)}(s)}(x))\,\mathrm{d}x,\,s\in[0,t], \label{sbp4}
\end{eqnarray} 
where $B_n$, $\Sigma_n$ are defined according to \eqref{BnSn}. Inserting the identities \eqref{sbp1}-\eqref{sbp4} into \eqref{eq_expansion1_26} we arrive at 
\begin{equation}\label{PMEhprelim}
\begin{split}
& \int_\rr \gamma(t,x)\,F_{\rho^{(n)}(t)}(x)\,\mathrm{d}x 
- \int_\rr \gamma(0,x)\,F_{\rho^{(n)}(0)}(x)\,\mathrm{d}x \\
& =-\frac{1}{n}\,\sum_{i=1}^n \int_0^t \gamma(s,X^{(n)}_i(s))\,\sigma(F_{\rho^{(n)}(s)}(X^{(n)}_i(s)))\,\mathrm{d}B^{(n)}_i(s) \\
&\quad\,+\!\int_0^t\!\int_\rr \big(\gamma_s(s,x) F_{\rho^{(n)}(s)}(x)\!+\!\gamma_x(s,x) B_n(F_{\rho^{(n)}(s)}(x))\!+\!\gamma_{xx}(s,x) \Sigma_n(F_{\rho^{(n)}(s)}(x))\big)\,\mathrm{d}x\,\mathrm{d}s.
\end{split}
\end{equation}

\smallskip

At this point, we take the difference between the equations \eqref{PMEhprelim} and \eqref{PMEh}, multiply the resulting equation by $\sqrt{n}$, use the Fundamental Theorem of Calculus in the forms 
\begin{eqnarray}
&& \begin{split}
B_n(F_{\rho^{(n)}(s)}(x))-B(R(s,x))
=B_n(F_{\rho^{(n)}(s)}(x))-B(F_{\rho^{(n)}(s)}(x))
\qquad\qquad\qquad\qquad \\ 
+\int_0^1 b(a F_{\rho^{(n)}(s)}(x)+(1-a) R(s,x))\,(F_{\rho^{(n)}(s)}(x)-R(s,x))\,\mathrm{d}a, 
\end{split} \\
&& \begin{split}
\Sigma_n(F_{\rho^{(n)}(s)}(x))-\Sigma(R(s,x))
=\Sigma_n(F_{\rho^{(n)}(s)}(x))-\Sigma(F_{\rho^{(n)}(s)}(x))
\qquad\qquad\qquad\qquad \\
+\int_0^1 \frac{\sigma(a F_{\rho^{(n)}(s)}(x)+(1-a) R(s,x))^2}{2}\,(F_{\rho^{(n)}(s)}(x)-R(s,x))\,\mathrm{d}a 
\end{split}
\end{eqnarray}
and rearrange terms to end up with \eqref{eq:main_repr}. \ep

\medskip

We are now ready to give the proof of Proposition \ref{prop:Gaussians}.

\medskip

\noindent\textbf{Proof of Proposition \ref{prop:Gaussians}. Step 1.} By definition the random variables in \eqref{limitGaussians} are the limits in law of 
\begin{equation}\label{Gaussians_prelimit}
\begin{split}
&\int_\rr \gamma(t_\ell,x)\,G_n(t_\ell)(\mathrm{d}x)
-\int_\rr \gamma(0,x)\,G_n(0)(\mathrm{d}x) \\
& -\int_0^{t_\ell} \int_\rr \Big(\gamma_s(s,x)+\gamma_x(s,x)\,b(R(s,x))+\gamma_{xx}(s,x)\,\frac{\sigma(R(s,x))^2}{2}\Big)\,H_n(t_\ell)(\mathrm{d}s,\mathrm{d}x), \\
&\int_\rr \gamma(0,x)\,G_n(0)(\mathrm{d}x) 
\end{split}
\end{equation}
along a suitable sequence of $n\in\nn$. 

\medskip

To proceed we note that the convergence $\rho^{(n)}\to\rho$ in probability in $C([0,\infty),M_1(\rr))$ and the regularity result of Proposition \ref{PME_grad} imply the convergences in probability 
\begin{equation}\label{unif_conv_probab}
\sup_{(s,x)\in[0,t_\ell]\times\rr} \big|F_{\rho^{(n)}(s)}(x)-R(s,x)\big|\to 0,\quad\ell=1,\,2,\,\ldots,\,m.
\end{equation}
This can be seen most easily by combining the Skorokhod Representation Theorem in the form of \cite[Theorem 3.5.1]{Du} for the sequence $\rho^{(n)}$, $n\in\nn$ with the regularity result of Proposition \ref{PME_grad} to first obtain the almost sure pointwise convergence $F_{\rho^{(n)}(s)}(\cdot)\to R(s,\cdot)$ for all $s\ge0$. Since all functions involved are cumulative distribution functions and one has the almost sure convergence $\rho^{(n)}\to\rho$ in $C([0,\infty),M_1(\rr))$, the almost sure convergence $F_{\rho^{(n)}(\cdot)}(\cdot)\to R(\cdot,\cdot)$ is in fact uniform on all sets of the form $[0,t]\times\rr$.

\medskip

The convergences of \eqref{unif_conv_probab} in conjunction with the Lipschitz property of $b$, $\frac{\sigma^2}{2}$ (cf. Assumption \ref{main_asmp}(b)) show that the limit in law of the random variables in \eqref{Gaussians_prelimit} along a sequence of $n\in\nn$ is the same as the limit in law of  
\begin{equation*}
\begin{split}
& \int_\rr \gamma(t_\ell,x)\,G_n(t_\ell)(\mathrm{d}x)-\int_\rr \gamma(0,x)\,G_n(0)(\mathrm{d}x) \\
&-\int_0^{t_\ell} \int_\rr \int_0^1 \bigg( 
\gamma_s(s,x)+\gamma_x(s,x)\,b(a F_{\rho^{(n)}(s)}(x)+(1-a)R(s,x)) \\
&\qquad\qquad\qquad\;\;
+\gamma_{xx}(s,x)\,\frac{\sigma(a F_{\rho^{(n)}(s)}(x)+(1-a) R(s,x))^2}{2}\bigg) \mathrm{d}a\,H_n(t_\ell)(\mathrm{d}s,\mathrm{d}x), \\
&\int_\rr \gamma(0,x)\,G_n(0)(\mathrm{d}x)
\end{split}
\end{equation*}
along the same sequence of $n\in\nn$.

\medskip

Next, we apply Lemma \ref{lemma:testfunction} and find that the latter limit in law must be equal to the limit in law of 
\begin{equation*}
\begin{split}
&-\frac{1}{\sqrt{n}}\,\sum_{i=1}^n \int_0^{t_\ell} \gamma(s,X^{(n)}_i(s))\,\sigma(F_{\rho^{(n)}(s)}(X^{(n)}_i(s)))\,\mathrm{d}B^{(n)}_i(s)\\
&+\sqrt{n}\int_0^{t_\ell} \int_\rr \big(\gamma_x(s,x)(B_n-B)(F_{\rho^{(n)}(s)}(x))+\gamma_{xx}(s,x)(\Sigma_n-\Sigma)(F_{\rho^{(n)}(s)}(x))\big)\,\mathrm{d}x\,\mathrm{d}s, \\
&\int_\rr \gamma(0,x)\,G_n(0)(\mathrm{d}x)
\end{split}
\end{equation*}
along the same sequence of $n\in\nn$. Moreover, since the functions $b$ and $\frac{\sigma^2}{2}$ are Lipschitz by Assumption \ref{main_asmp}(b), the suprema $\sup_{[0,1]} |B_n-B|$ and $\sup_{[0,1]} |\Sigma_n-\Sigma|$ can be bounded above by $C n^{-1}$ with a constant $C<\infty$ depending only on the Lipschitz constants of $b$ and $\frac{\sigma^2}{2}$. Consequently, it suffices to study the limit in law of 
\begin{equation}
-\frac{1}{\sqrt{n}}\,\sum_{i=1}^n \int_0^{t_\ell} \gamma(s,X^{(n)}_i(s))\,\sigma(F_{\rho^{(n)}(s)}(X^{(n)}_i(s)))\,\mathrm{d}B^{(n)}_i(s),\quad
\int_\rr \gamma(0,x)\,G_n(0)(\mathrm{d}x)
\end{equation}
along the same sequence of $n\in\nn$ as before.

\medskip

\noindent\textbf{Step 2.} Consider the sequences of continuous martingales
\begin{equation}\label{cont_mart_def}
\int_\rr \gamma(0,x)\,G_n(0)(\mathrm{d}x)-\frac{1}{\sqrt{n}}\,\sum_{i=1}^n \int_0^t \gamma(s,X^{(n)}_i(s))\,\sigma(F_{\rho^{(n)}(s)}(X^{(n)}_i(s)))\,\mathrm{d}B^{(n)}_i(s),\,t\in[0,t_\ell]
\end{equation}
indexed by $n\in\nn$, where $\ell$ and $\gamma$ vary over $\{1,2,\ldots,m\}$ and a countable dense subset ${\mathcal C}_\ell$ of the space of functions on $[0,t_\ell]\times\rr$ which are continuously differentiable in $s$, twice continuously differentiable in $x$ and compactly supported. One easily verifies the tightness of each such sequence via the tightness criterion of \cite[Theorem 7.3]{Bi} by recalling Proposition \ref{CLT_dGM}, writing each of the martingales as a time-changed standard Brownian motion with the same initial value (cf. \cite[Chapter 3, Problem 4.7]{KS}) and using the assumed boundedness of $\gamma$ and $\sigma$. In particular, every sequence of $n\in\nn$ admits a subsequence along which the continuous martingales of \eqref{cont_mart_def} converge to the respective limiting processes $M^\gamma$ for all $\gamma\in{\mathcal C}_\ell$, $\ell\in\{1,2,\ldots,m\}$. 

\medskip

Now, letting $\Gamma(s,x):=-\int_x^\infty \gamma(s,y)\,\mathrm{d}y$ as before, integrating by parts, recalling Assumption \ref{main_asmp}(a), applying the inequality \eqref{ineq_elem_p} with $p=2$ and using the It\^o isometry we arrive at the estimate
\begin{equation}
\begin{split}
& 2\,\Big(\E\big[\Gamma(0,X^{(n)}_1(0))^2\big]
-\E\big[\Gamma(0,X^{(n)}_1(0))\big]^2\Big) \\
& +2\,\E\bigg[\int_0^t \int_\rr \gamma(s,x)^2\, \sigma(F_{\rho^{(n)}(s)}(x))^2\,\rho^{(n)}(s)(\mathrm{d}x)\,\mathrm{d}s\bigg]
\end{split}
\end{equation}
on the second moment of the random variable in \eqref{cont_mart_def} with the same value of $t$. The latter quantities tend to 
\begin{equation}\label{eq:QV_lim}
2\,\Big(\E\big[\Gamma(0,X^{(1)}_1(0))^2\big]-\E\big[\Gamma(0,X^{(1)}_1(0))\big]^2\Big)
+2\,\E\bigg[\int_0^t \int_\rr \gamma(s,x)^2\,\sigma(R(s,x))^2\,\rho(s)(\mathrm{d}x)\,\mathrm{d}s\bigg]
\end{equation}
in the limit $n\to\infty$, as can be seen by applying the Skorokhod Embedding Theorem in the form of \cite[Theorem 3.5.1]{Du} to the sequence $\rho^{(n)}$, $n\in\nn$, using the almost sure weak convergences 
\begin{equation}\label{aswc}
\begin{split}
\sigma(F_{\rho^{(n)}(s)}(x))^2\,\rho^{(n)}(s)(\mathrm{d}x)\!=\!
2\,\mathrm{d}\Sigma_n(F_{\rho^{(n)}(s)}(\cdot))\!\to\! 
2\,\mathrm{d}\Sigma(R(s,\cdot))
\!=\!\sigma(R(s,x))^2\,\rho(s)(\mathrm{d}x), \\
s\in[0,t]
\end{split}
\end{equation}
and appealing to the Dominated Convergence Theorem (recall that $\gamma$ and $\sigma$ are bounded by assumption). In particular, the one-dimensional distributions of the continuous martingales in \eqref{cont_mart_def} are uniformly integrable as $n$ varies, so that the limiting processes $M^\gamma$ must be themselves continuous martingales for all $\gamma\in{\mathcal C}_\ell$, $\ell\in\{1,2,\ldots,m\}$. 

\medskip

Finally, for any $\gamma\in{\mathcal C}_\ell$, $\tilde{\gamma}\in{\mathcal C}_{\tilde{\ell}}$ another application of the Skorokhod Embedding Theorem to the sequence $\rho^{(n)}$, $n\in\nn$, the convergences in \eqref{aswc} and the Dominated Convergence Theorem show that the quadratic covariation process on $[0,\min(t_\ell,t_{\tilde{\ell}})]$ between the continuous martingales of \eqref{cont_mart_def} associated with $\gamma$, $\tilde{\gamma}$ converges in law to 
\begin{equation}\label{covar_limit}
\int_0^t \int_\rr \gamma(s,x)\,\tilde{\gamma}(s,x)\,\sigma(R(s,x))^2\,\rho(s)(\mathrm{d}x)\,\mathrm{d}s,\quad t\in[0,\min(t_\ell,t_{\tilde{\ell}})]
\end{equation}
in the limit $n\to\infty$. Moreover, another uniform integrability argument relying on integration by parts, Assumption \ref{main_asmp}(a), the inequality \eqref{ineq_elem_p} with $p=4$,  the Burkholder-Davis-Gundy inequality (see e.g. \cite[Chapter 3, Theorem 3.28]{KS}) and the boundedness of $\gamma$ and $\sigma$ allows to identify the process in \eqref{covar_limit} as the quadratic covariation process between $M^\gamma$ and $M^{\tilde{\gamma}}$. This and Proposition \ref{CLT_dGM} lead to the conclusion that the probability space supporting $M^\gamma$, $\gamma\in{\mathcal C}_\ell$, $\ell\in\{1,2,\ldots,m\}$ admits an orthogonal martingale measure ${\mathrm d}M(s,x)$ on $[0,t_m]\times\rr$ in the sense of \cite[definitions on pp. 287--288]{Wa} with the quadratic variation measure 
\begin{equation}
\mathrm{d}\langle M\rangle(s,x)=\sigma(R(s,x))^2\,\rho(s)(\mathrm{d}x)\,\mathrm{d}s\quad\text{on}\quad [0,t_m]\times\rr
\end{equation}
and a reparametrized Brownian bridge $\beta(F_\lambda(\cdot))$ independent of ${\mathrm d}M(s,x)$ satisfying
\begin{equation}
M^\gamma(t)=\int_\rr \gamma(0,x)\,\beta(F_\lambda(x))\,\mathrm{d}x+\int_0^t \int_\rr \gamma(s,x)\,\mathrm{d}M(s,x),\quad t\in[0,t_\ell]
\end{equation}
for all $\gamma\in{\mathcal C}_\ell$, $\ell\in\{1,2,\ldots,m\}$. It remains to use the positivity of $\sigma$ throughout $[0,1]$ and the existence of a positive density $R_x(s,\cdot)$ of $\rho(s)$ for $s>0$ (cf. \eqref{lawbarX} and the lower bound of \eqref{eq:heat_kernel}) in order to define the white noise 
\begin{equation}\label{white_noise_coupling}
\mathrm{d}W(s,x):=\sigma(R(s,x))^{-1}\,R_x(s,x)^{-1/2}\,\mathrm{d}M(s,x)\quad\text{on}\quad [0,t_m]\times\rr,
\end{equation} 
ending up with the identification 
\begin{equation}
M^\gamma(t)\!=\!\int_\rr \gamma(0,x)\,\beta(F_\lambda(x))\,\mathrm{d}x+\int_0^t \! \int_\rr \gamma(s,x)\,\sigma(R(s,x))\,R_x(s,x)^{1/2}\,\mathrm{d}W(s,x),\; t\in[0,t_\ell]
\end{equation}
for all $\gamma\in{\mathcal C}_\ell$, $\ell\in\{1,2,\ldots,m\}$. The statement of the proposition for such $\ell$ and $\gamma$ readily follows. To obtain the statement for arbitrary $\ell$ and $\gamma$ it suffices to pick a sequence of functions from ${\mathcal C}_\ell$ converging to $\gamma$, use the statement for the latter and pass to the limit. \ep

\medskip

We proceed to an analogue of Proposition \ref{prop:Gaussians} for the mild solution $G$ from \eqref{mild_def}.

\begin{proposition}\label{prop:mild_Gaussians}
Suppose that Assumption \ref{main_asmp} holds. Then, for any $t>0$ the measures $G(t,x)\,\mathrm{d}x$ on $\rr$ and $G(s,x)\,\mathbf{1}_{[0,t]\times\rr}(s,x)\,\mathrm{d}s\,\mathrm{d}x$ on $[0,t]\times\rr$, defined in terms of the mild solution $G$ from \eqref{mild_def}, are finite almost surely and for every function on $[0,t]\times\rr$ which is continuously differentiable in $s$, twice continuously differentiable in $x$ and compactly supported one has
\begin{equation}\label{eq:mild_Gaussians}
\begin{split}
& \int_\rr \gamma(t,x)\,G(t,x)\,\mathrm{d}x
-\int_\rr \gamma(0,x)\,G(0,x)\,\mathrm{d}x \\
& -\int_0^t \int_\rr \Big(\gamma_s(s,x)+\gamma_x(s,x)\,b(R(s,x))+\gamma_{xx}(s,x)\,\frac{\sigma(R(s,x))^2}{2}\Big)\,G(s,x)\,\mathrm{d}x\,\mathrm{d}s\\
& =\int_0^t \int_\rr \gamma(s,x)\,\sigma(R(s,x))\,R_x(s,x)^{1/2}\,\mathrm{d}W(s,x).
\end{split}
\end{equation}
\end{proposition}

\medskip

\noindent\textbf{Proof. Step 1.} We fix a $t>0$ and aim to verify in this first step that
\begin{equation}\label{mG0}
\E\bigg[\int_\rr |G(t,x)|\,\mathrm{d}x\bigg]<\infty\quad\text{and}\quad \E\bigg[\int_0^t \int_\rr |G(s,x)|\,\mathrm{d}x\,\mathrm{d}s\bigg]<\infty.
\end{equation}
To this end, we insert the right-hand side of \eqref{mild_def} into the first expectation and bound the result using the triangle inequality, Fubini's Theorem and Jensen's inequality by
\begin{equation}\label{mG1}
\E\bigg[\int_\rr \big|\beta(F_\lambda(y))\big|\,\mathrm{d}y\bigg]
+\int_\rr \bigg(\int_0^t \int_\rr \sigma(R(s,y))^2\,R_x(s,y)\,p(s,y;t,x)^2\,\mathrm{d}y\,\mathrm{d}s\bigg)^{1/2}\,\mathrm{d}x.
\end{equation}
Fubini's Theorem and the scaling property of Gaussian distributions reveals further that the first summand in \eqref{mG1} is the product of the first absolute moment of the standard Gaussian distribution and $\int_\rr \sqrt{F_\lambda(y)(1-F_\lambda(y))}\,\mathrm{d}y$. The latter integral is finite due to Assumption \ref{main_asmp}(a) and \cite[discussion of the functional $J_1$ on p. 25]{BL}.

\medskip

To estimate the second summand in \eqref{mG1} we combine the boundedness of $\sigma$ (cf. Assumption \ref{main_asmp}(b)), the inequality $p(s,y;t,x)\le C(t-s)^{-1/2}$ (cf. \eqref{eq:heat_kernel}) and the identity
\begin{equation}
\int_\rr R_x(s,y)\,p(s,y;t,x)\,\mathrm{d}y=R_x(t,x)
\end{equation}
(due to the Markov property of the diffusion $\bar{X}$, see Subsection \ref{sec:PMEetc}) to arrive at the upper bound 
\begin{equation}
C\,\int_\rr \bigg(\int_0^t \frac{1}{2}\,(t-s)^{-1/2}\,R_x(t,x)\,\mathrm{d}s\bigg)^{1/2}\,\mathrm{d}x=C\,t^{1/4}\,\int_\rr R_x(t,x)^{1/2}\,\mathrm{d}x,
\end{equation}
where $C<\infty$ depends only on $\sup_{[0,1]} \sigma$ and the constant in \eqref{eq:heat_kernel}. At this point, Jensen's inequality with respect to the Cauchy distribution 
\begin{equation}
\int_\rr\!R_x(t,x)^{1/2}\,\mathrm{d}x
\!=\!\pi\!\int_\rr\! R_x(t,x)^{1/2}(1\!+\!x^2)\frac{1}{\pi(1\!+\!x^2)}\mathrm{d}x
\!\le\!\pi^{1/2}\bigg(\!\int_\rr\! R_x(t,x)(1\!+\!x^2)\,\mathrm{d}x\!\bigg)^{\!1/2}
\end{equation}
and the estimate \eqref{moment_est} imply that the first expectation in \eqref{mG0} is finite. Moreover, in view of Fubini's Theorem and since the just obtained estimate is uniformly bounded on every compact interval of $t$'s, the second expectation in \eqref{mG0} is also finite. 
 
\medskip

\noindent\textbf{Step 2.} To derive the identity \eqref{eq:mild_Gaussians} we fix a function $\gamma$ as described and deduce from the definition of $G$ in \eqref{mild_def} that
\begin{equation}\label{mG2}
\begin{split}
&\int_\rr \gamma(t,x)\,G(t,x)\,\mathrm{d}x 
= \int_\rr \gamma(t,x) \int_\rr G(0,y)\,p(0,y;t,x)\,\mathrm{d}y\,\mathrm{d}x \\
&\qquad\qquad\qquad\qquad\quad
+\int_\rr \gamma(t,x) \int_0^t\int_\rr \sigma(R(s,y))\,R_x(s,y)^{1/2}\,p(s,y;t,x)\,\mathrm{d}W(s,y)\,\mathrm{d}x.
\end{split}
\end{equation}
Moreover, the boundedness of $\gamma$ and $\sigma$ and the estimates
\begin{eqnarray}
&&\quad\int_\rr \int_\rr |G(0,y)|\,p(0,y;t,x)\,\mathrm{d}y\,\mathrm{d}x
=\int_\rr |G(0,y)|\,\mathrm{d}y<\infty, \\
&&\quad\int_\rr \int_0^t \int_\rr R_x(s,y)\,p(s,y;t,x)^2\,\mathrm{d}y\,\mathrm{d}s\,\mathrm{d}x\le C\int_\rr \int_0^t (t-s)^{-1/2}
R_x(t,x)\,\mathrm{d}s\,\mathrm{d}x<\infty
\end{eqnarray}
(see Step 1 for more details) allow us to use the classical and the stochastic Fubini's Theorems (see \cite[Theorem 2.6]{Wa} and note that the dominating measure therein is $\delta_{\tilde{y}}(\mathrm{d}y)\,\mathrm{d}\tilde{y}\,\mathrm{d}s$ in our case) and to rewrite the right-hand side of \eqref{mG2} as 
\begin{equation}\label{mG3}
\begin{split}
& \int_\rr G(0,y) \int_\rr \gamma(t,x)\,p(0,y;t,x)\,\mathrm{d}x\,\mathrm{d}y \\
& +\int_0^t\int_\rr \sigma(R(s,y))\,R_x(s,y)^{1/2}\,\int_\rr \gamma(t,x)\,p(s,y;t,x)\,\mathrm{d}x\,\mathrm{d}W(s,y).
\end{split}
\end{equation}

\smallskip

Next, we employ It\^o's formula and Fubini's Theorem to find
\begin{equation*}
\begin{split}
\int_\rr \gamma(t,x)\,p(s,y;t,x)\,\mathrm{d}x &=
\E\big[\gamma(t,\bar{X}(t))\big|\bar{X}(s)=y\big] \\
&=\gamma(s,y)+\E\bigg[\int_s^t ({\mathcal A}_r\gamma)(r,\bar{X}(r))\,\mathrm{d}r\bigg|\bar{X}(s)=y\bigg] \\
&=\gamma(s,y)+\int_s^t \int_\rr ({\mathcal A}_r\gamma)(r,x)\,p(s,y;r,x)\,\mathrm{d}x\,\mathrm{d}r,
\end{split}
\end{equation*} 
where 
\begin{equation}\label{Argamma}
({\mathcal A}_r\gamma)(r,x):=\gamma_s(r,x)+\gamma_x(r,x)\,b(R(r,x))+\gamma_{xx}(r,x)\,\frac{\sigma(R(r,x))^2}{2},\;\; (r,x)\in[0,t]\times\rr.
\end{equation}
Applying this observation to the expression in \eqref{mG3} we get
\begin{equation}\label{mG4}
\begin{split}
&\int_\rr G(0,y)\,\gamma(0,y)\,\mathrm{d}y
+\int_\rr G(0,y) \int_0^t \int_\rr  ({\mathcal A}_r \gamma)(r,x)\, p(0,y;r,x)\,\mathrm{d}x\,\mathrm{d}r\,\mathrm{d}y \\
&+\int_0^t \int_\rr \sigma(R(s,y))\,R_x(s,y)^{1/2}\,\gamma(s,y)\,\mathrm{d}W(s,y) \\
&+\int_0^t \int_\rr \sigma(R(s,y))\,R_x(s,y)^{1/2} \int_s^t
\int_\rr ({\mathcal A}_r \gamma)(r,x)\,p(s,y;r,x)\,\mathrm{d}x\,\mathrm{d}r\,\mathrm{d}W(s,y).
\end{split}
\end{equation}

\smallskip

At this point, thanks to the boundedness of $({\mathcal A}_r\gamma)$ and $\sigma$ (cf. Assumption \ref{main_asmp}(b)) and the estimates
\begin{eqnarray}
&&\int_\rr |G(0,y)| \int_0^t \int_\rr p(0,y;r,x)\,\mathrm{d}x\,\mathrm{d}r\,\mathrm{d}y=\int_\rr |G(0,y)|\,t\,\mathrm{d}y<\infty, \\
&&\begin{split}
& \int_0^t \int_\rr \int_0^r \int_\rr R_x(s,y)\,p(s,y;r,x)^2\,\mathrm{d}y\,\mathrm{d}s\,\mathrm{d}x\,\mathrm{d}r \\
& \le C\int_0^t\int_\rr\int_0^r (r-s)^{-1/2}\,R_x(r,x)\,\mathrm{d}s\,\mathrm{d}x\,\mathrm{d}r<\infty  
\end{split}
\end{eqnarray}
the classical and the stochastic Fubini's Theorems are applicable to the second and fourth summands in \eqref{mG4}, so that the overall expression in \eqref{mG4} equals to
\begin{equation}
\begin{split}
&\int_\rr \gamma(0,y)\,G(0,y)\,\mathrm{d}y
+\int_0^t \int_\rr ({\mathcal A}_r \gamma)(r,x) \int_\rr
G(0,y)\,p(0,y;r,x)\,\mathrm{d}y\,\mathrm{d}x\,\mathrm{d}r \\
&+\int_0^t \int_\rr \sigma(R(s,y))\,R_x(s,y)^{1/2}\,\gamma(s,y)\,\mathrm{d}W(s,y) \\
&+\int_0^t \int_\rr ({\mathcal A}_r \gamma)(r,x)
\int_0^r \int_\rr \sigma(R(s,y))\,R_x(s,y)^{1/2}\,p(s,y;r,x)\,\mathrm{d}W(s,y)\,\mathrm{d}x\,\mathrm{d}r \\
&=\int_\rr \gamma(0,y)\,G(0,y)\,\mathrm{d}y
+\int_0^t \int_\rr ({\mathcal A}_r \gamma)(r,x)\,G(r,x)\,\mathrm{d}x\,\mathrm{d}r \\
&\quad+\int_0^t \int_\rr \gamma(s,y)\,\sigma(R(s,y))\,R_x(s,y)^{1/2}\,\mathrm{d}W(s,y).
\end{split}
\end{equation} 
This finishes the proof of the proposition. \ep 

\medskip

We can now identify the subsequential limits of Proposition \ref{prop:tight}.

\begin{proposition}\label{prop:full_iden}
Suppose that Assumption \ref{main_asmp} is satisfied. Then, any subsequential limit in law of the sequence in \eqref{FDDseq} has the same distribution as
\begin{equation}\label{eq:limitingrandomvector}
\begin{split}
& \big(G(0,x)\,\mathrm{d}x,\,G(t_1,x)\,\mathrm{d}x,\,\ldots,\, G(t_m,x)\,\mathrm{d}x, \\
&G(s,x)\mathbf{1}_{[0,t_1]\times\rr}(s,x)\,\mathrm{d}s\,\mathrm{d}x,\,G(s,x)\mathbf{1}_{[0,t_2]\times\rr}(s,x)\,\mathrm{d}s\,\mathrm{d}x,\,\ldots,\,G(s,x)\mathbf{1}_{[0,t_m]\times\rr}(s,x)\,\mathrm{d}s\,\mathrm{d}x\big),
\end{split}
\end{equation}
where $G$ is the mild solution from \eqref{mild_def}.
\end{proposition}

\noindent\textbf{Proof. Step 1.} We consider a probability space that supports a limit point in law 
\begin{equation}
\big(G_\infty(0), G_\infty(t_1), \ldots, G_\infty(t_m), H_\infty(t_1), H_\infty(t_2), \ldots, H_\infty(t_m)\big)
\end{equation}
of the sequence in \eqref{FDDseq} and aim to couple it with a mild solution of the SPDE \eqref{SPDE}. 

\medskip

To this end, for each $\ell\in\{1,2,\ldots,m\}$ we pick a countable dense subset ${\mathcal C}_\ell$ of the space of  
functions on $[0,t_\ell]\times\rr$ which are continuously differentiable in $s$, twice continuously differentiable in $x$ and compactly supported. We note that the random variables of \eqref{limitGaussians} with $\ell$ and $\gamma$ varying over $\{1,2,\ldots,m\}$ and ${\mathcal C}_\ell$ are defined on the underlying probability space. Moreover, by Proposition \ref{prop:Gaussians} their joint distribution must be that of the random variables in \eqref{limitGaussians'}. Hence, by \cite[Theorem 5.3]{Ka} we can define on an enlargement of the underlying probability space a countable collection of continuous processes whose conditional distribution given the random variables in \eqref{limitGaussians} with $\ell$ and $\gamma$ varying over $\{1,2,\ldots,m\}$ and ${\mathcal C}_\ell$ is the same as the conditional distribution of the continuous processes  
\begin{equation}
\int_0^t \int_\rr \gamma(s,x)\,\sigma(R(s,x))\,R_x(s,x)^{1/2}\,\mathrm{d}W(s,x),\;\;t\in[0,t_\ell],\quad \gamma\in{\mathcal C}_\ell,\;\;\ell=1,\,2,\,\ldots,\,m
\end{equation}
given
\begin{equation}
\begin{split}
& \int_0^{t_\ell} \int_\rr \gamma(s,x)\,\sigma(R(s,x))\,R_x(s,x)^{1/2} \,\mathrm{d}W(s,x),\;\;\gamma\in{\mathcal C}_\ell,\;\;\ell=1,\,2,\,\ldots,\,m, \\
& \int_\rr \gamma(0,x)\,\beta(F_\lambda(x))\,\mathrm{d}x,\;\;\gamma\in{\mathcal C}_\ell,\;\;\ell=1,\,2,\,\ldots,\,m.
\end{split}
\end{equation}

\smallskip

It follows that the enlarged probability space supports an orthogonal martingale measure ${\mathrm d}M(s,x)$ on $[0,t_m]\times\rr$ in the sense of \cite[definitions on pp. 287--288]{Wa} with the quadratic variation measure 
\begin{equation}
\mathrm{d}\langle M\rangle(s,x)=\sigma(R(s,x))^2\,R_x(s,x)\,\mathrm{d}x\,\mathrm{d}s\quad\text{on}\quad [0,t_m]\times\rr
\end{equation}
and we can define a white noise $\mathrm{d}W(s,x)$ on $[0,t_m]\times\rr$ as in \eqref{white_noise_coupling}. Finally, we let $G$ be the mild solution of the SPDE \eqref{SPDE} on $[0,t_m]\times\rr$ given by 
\begin{equation}
\begin{split}
G(t,x)=\int_\rr p(0,y;t,x)\,G_\infty(0)(\mathrm{d}y)
+\int_0^t \!\int_\rr \sigma(R(s,y))\,R_x(s,y)^{1/2}\,p(s,y;t,x)\,\mathrm{d}W(s,y), \\
(t,x)\in[0,t_m]\times\rr.
\end{split}
\end{equation}
In particular, Proposition \ref{prop:mild_Gaussians} and our coupling construction ensure that
\begin{equation}\label{eq:Delta}
\begin{split}
&\int_\rr \gamma(t_\ell,x)\,G(t_\ell,x)\,\mathrm{d}x\!-\!\int_0^{t_\ell}\!\int_\rr \,({\mathcal A}_s\gamma)(s,x)\,G(s,x)\,\mathrm{d}x\,\mathrm{d}s \\
&=\!\int_\rr \gamma(t_\ell,x)\,G_\infty(t_\ell)(\mathrm{d}x)
\!-\!\int_0^{t_\ell} \!\int_\rr \,({\mathcal A}_s\gamma)(s,x)\,H_\infty(t_\ell)(\mathrm{d}s,\mathrm{d}x),\;\gamma\in{\mathcal C}_\ell,\;\ell=1,\,2,\,\ldots,\,m,
\end{split}
\end{equation}
with the notation of \eqref{Argamma}.

\medskip

\noindent\textbf{Step 2.} We fix an $\ell\in\{1,2,\ldots,m\}$ and a continuous function $g:\,[0,t_\ell]\times\rr\to\rr$ with compact support and consider the backward Cauchy problem
\begin{equation}\label{Krylov_PDE2}
{\mathcal A}_s u=g,\quad u(t_\ell,\cdot)=0
\end{equation} 
on $[0,t_\ell]\times\rr$. As explained in the paragraph following \eqref{Krylov_PDE}, the conditions of \cite[Theorem 2.1]{Kr2} apply to the equation \eqref{Krylov_PDE2} and guarantee the existence of a solution $u$ with $u,u_t,u_x,u_{xx}\in L^2([0,t_\ell]\times\rr)$. We claim that \eqref{eq:Delta} implies
\begin{equation}\label{eq:Delta2}
\begin{split}
&\int_\rr u(t_\ell,x)\,G(t_\ell,x)\,\mathrm{d}x-\int_0^{t_\ell} \int_\rr \,({\mathcal A}_s u)(s,x)\,G(s,x)\,\mathrm{d}x\,\mathrm{d}s \\
&=\int_\rr u(t_\ell,x)\,G_\infty(t_\ell)(\mathrm{d}x)
-\int_0^{t_\ell} \int_\rr \,({\mathcal A}_su)(s,x)\,H_\infty(t_\ell)(\mathrm{d}s,\mathrm{d}x).
\end{split}
\end{equation}  
Since the first integrals on both sides of \eqref{eq:Delta2} vanish due to the terminal condition in \eqref{Krylov_PDE2} and $g={\mathcal A}_s u$ can be chosen arbitrarily from a countable dense subset of $C_0([0,t_\ell]\times\rr)$, it would follow from \eqref{eq:Delta2} that $G(s,x)\,\mathbf{1}_{[0,t_\ell]\times\rr}\,\mathrm{d}s\,\mathrm{d}x=H_\infty(t_\ell)(\mathrm{d}s,\mathrm{d}x)$ for all $\ell\in\{1,2,\ldots,m\}$ and then from \eqref{eq:Delta} that $G(t_\ell,x)\,\mathrm{d}x=G_\infty(t_\ell)(\mathrm{d}x)$ for all $\ell\in\{1,2,\ldots,m\}$, finishing the proof of the proposition.

\medskip 

To obtain \eqref{eq:Delta2} from \eqref{eq:Delta} it suffices to show that one can pick functions $\gamma^{(\kappa)}$, $\kappa\in\nn$ in ${\mathcal C}_\ell$ with
\begin{equation}\label{approx_by_smooth}
\gamma^{(\kappa)}(t_\ell,\cdot)\to u(t_\ell,\cdot)\;\;\text{and}\;\;{\mathcal A}_s \gamma^{(\kappa)}\to {\mathcal A}_s u=g\;\;\text{uniformly as}\;\;\kappa\to\infty.
\end{equation}
To this end, we recall the solution $\bar{X}$ of the SDE \eqref{PME_SDE} and observe that the time-homogeneous Markov process $(s,\bar{X}(s))$, $s\in[0,t_\ell]$ is the unique weak solution of the associated SDE, hence also of the local martingale problem for the operator ${\mathcal A}_s$ (see e.g. \cite[Theorem 18.7]{Ka}). The latter has bounded continuous coefficients, so that $(s,\bar{X}(s))$, $s\in[0,t_\ell]$ is a Feller process and its generator is the unique extension of ${\mathcal A}_s$ from the space of infinitely differentiable functions with compact support in $[0,t_\ell]\times\rr$ to an appropriate domain within the space of continuous functions on $[0,t_\ell]\times\rr$ vanishing at infinity (see e.g. \cite[Theorem 18.11]{Ka}).  

\medskip

Next, we employ the stochastic representation 
\begin{equation}\label{stoch_repr}
u(s,x)=-\int_s^t \E\big[g(r,\bar{X}(r))\big|\bar{X}(s)=x\big]\,\mathrm{d}r,\quad (s,x)\in[0,t_\ell]\times\rr
\end{equation} 
of the solution to \eqref{Krylov_PDE2} (cf. the explanation preceding \eqref{Feynman-Kac}). Together with the Feller property of the process $(s,\bar{X}(s))$, $s\in[0,t]$ and the Dominated Convergence Theorem it shows that $u$ is continuous. Moreover, since $g$ has compact support and the diffusion $\bar{X}$ has bounded coefficients, $u$ vanishes at infinity. Finally, the representation \eqref{stoch_repr} reveals that the process
\begin{equation}
u(s,\bar{X}(s))-u(0,\bar{X}(0))-\int_0^s g(r,\bar{X}(r))\,\mathrm{d}r,\quad s\in[0,t]
\end{equation}
is a martingale, so that by the converse of Dynkin's formula (see e.g. \cite[Chapter VII, Proposition 1.7]{RY}) $u$ belongs to the domain of ${\mathcal A}_s$ with ${\mathcal A}_su=g$. In particular, $u$ admits an approximation as described in \eqref{approx_by_smooth}. \ep 

\medskip

We conclude the section with the proof of Theorem \ref{thm:CLT}.

\medskip

\noindent\textbf{Proof of Theorem \ref{thm:CLT}.} By Proposition \ref{prop:tight} every subsequence of the sequence in \eqref{FDDseq} has a further subsequence which converges in law. Moreover, by Proposition \ref{prop:full_iden} the limit of the latter must have the distribution of the random vector in \eqref{eq:limitingrandomvector}. Consequently, the whole sequence in \eqref{FDDseq} converges in law to the random vector in \eqref{eq:limitingrandomvector}, which is precisely the content of Theorem \ref{thm:CLT}. \ep

\bigskip\bigskip\bigskip



\begin{thebibliography}{BFPS}

\bibitem[Ar]{Ar} D.~G.~Aronson (1968). Non-negative solutions of linear parabolic equations. \textit{Ann. Scuola Norm. Sup. Pisa} \textbf{22}, pp. 607-–694.

\bibitem[BP]{BP} R.~F.~Bass, E.~Pardoux (1987). Uniqueness for diffusions with piecewise constant coefficients. \textit{Probab. Theory Related Fields} \textbf{76}, pp. 557--572. 

\bibitem[Bi]{Bi} P.~Billingsley (1999). \textit{Convergence of probability measures}. 2nd ed. John Wiley \& Sons, New York. 

\bibitem[BL]{BL} S.~Bobkov, M.~Ledoux (2014). One-dimensional empirical measures, order statistics and Kantorovich transport distances. Preprint available at \textit{math.umn.edu/$\sim$bobko001/preprints/2014\_BL\_Order.statistics.13.pdf}.

\bibitem[CP]{CP} S.~Chatterjee, S.~Pal (2011). A combinatorial analysis of interacting diffusions. \textit{J. Theoret. Probab.} \textbf{24}, pp. 939--968.

\bibitem[DG]{DG} D.~A.~Dawson, J.~G\"artner (1987). Large deviations from the McKean-Vlasov limit for weakly interacting diffusions. \textit{Stochastics} \textbf{20}, pp. 247--308.

\bibitem[dGM]{dGM} E.~del Barrio, E.~Gin\'{e}, C.~Matr\'{a}n (1999). Central limit theorems for the Wasserstein distance between the empirical and the true distributions. \textit{Ann. Probab.} \textbf{27}, pp. 1009--1071. 

\bibitem[DSVZ]{DSVZ} A.~Dembo, M.~Shkolnikov, S.~R.~S.~Varadhan, O.~Zeitouni (2016). Large deviations for diffusions interacting through their ranks. \textit{Comm. Pure Appl. Math.} \textbf{69}, pp. 1259--1313.

\bibitem[DK]{DK} J.~I.~Diaz, R.~Kersner (1987). On a nonlinear degenerate parabolic equation in infiltration or evaporation
through a porous medium. \textit{J. Differential Equations} \textbf{69}, pp. 368-–403.

\bibitem[Du]{Du} R.~M.~Dudley (1999). \textit{Uniform central limit theorems}. Cambridge University Press.

\bibitem[Fe]{Fe} E.~R.~Fernholz (2002). \textit{Stochastic portfolio theory}. Applications of Mathematics \textbf{48}. Springer-Verlag, New York.

\bibitem[FK]{FK} R.~Fernholz, I.~Karatzas (2009). Stochastic portfolio theory: an overview. In: A.~Bensoussan, Q.~Zhang (eds.) \textit{Handbook of Numerical Analysis}. Mathematical
Modeling and Numerical Methods in Finance \textbf{XV}, pp. 89--167. North-Holland, Oxford.

\bibitem[FGH]{FGH} D.~H.~Fremlin, D.~J.~H.~Garling, R.~G.~Haydon (1972). Bounded measures on topological spaces. \textit{Proc. London Math. Soc.} \textbf{25}, pp. 115--136.

\bibitem[Ga]{Ga} J.~G\"artner (1988). On the McKean-Vlasov limit for interacting diffusions. \textit{Math. Nachr.} \textbf{137}, pp. 197--248.

\bibitem[Gi]{Gi} B.~H.~Gilding (1989). Improved theory for a nonlinear degenerate parabolic equation. \textit{Ann. Sc. Norm. Super. Pisa Cl. Sci.} \textbf{16}, pp. 165--224.

\bibitem[IPS]{IPS} T.~Ichiba, S.~Pal, M.~Shkolnikov (2013). Convergence rates for rank-based models with applications to portfolio theory. \textit{Probab. Theory Related Fields} \textbf{156}, pp. 415--448.

\bibitem[JR]{JR} B.~Jourdain, J. Reygner (2013). Propagation of chaos for rank-based interacting diffusions and long time behaviour of a scalar quasilinear parabolic equation. \textit{Stochastic Partial Differential Equations: Analysis and Computations} \textbf{1}, pp. 455--506.

\bibitem[Ka]{Ka} O.~Kallenberg (1997). \textit{Foundations of Modern Probability}. Springer-Verlag, New York. 

\bibitem[KS]{KS} I.~Karatzas, S.~Shreve (1991). \textit{Brownian motion and stochastic calculus}. 2nd ed. Springer-Verlag, New York.

\bibitem[Kr1]{Kr1} N.~V.~Krylov (1971). An inequality in the theory of stochastic integrals. \textit{Theory Probab. Appl.} \textbf{16}, pp. 438--448. 

\bibitem[Kr2]{Kr2} N.~V.~Krylov (2007). Parabolic and elliptic equations with VMO coefficients. \textit{Comm. Partial Differential Equations} \textbf{32}, pp. 453--475.

\bibitem[La]{La} P.~D.~Lax (2002). \textit{Functional analysis}.  John Wiley \& Sons, New York.


\bibitem[Le]{Le} C.~Leonard (1986). Une loi des grands nombres pour des syst\`emes de diffusions avec interaction et \`a coefficients non born\'{e}s. \textit{Ann. Inst. Henri Poincar\'{e}} \textbf{22}, pp. 237--262



\bibitem[Mc]{Mc} H.~P.~McKean (1969). Propagation of chaos for a class of nonlinear parabolic equations. In: Lecture Series in Differential Equations \textbf{2}, pp. 41--57. Van Nostrand. Reinhold Co., New York. 

\bibitem[Oe1]{Oe1} K.~Oelschl\"ager (1984). A martingale approach to the law of large numbers for weakly interacting stochastic processes. \textit{Ann. Probab.} \textbf{12}, pp. 458--479.

\bibitem[Oe2]{Oe2} K.~Oelschl\"ager (1987). A fluctuation theorem for moderately interacting diffusion processes. \textit{Probab. Th. Rel. Fields} \textbf{74}, pp. 591--616. 

\bibitem[RY]{RY} D.~Revuz, M.~Yor (1999). \textit{Continuous martingales and Brownian motion}. Springer-Verlag, New York.

\bibitem[S]{S} M.~Shkolnikov (2012). Large systems of diffusions interacting through their ranks. \textit{Stochastic Process. Appl.} \textbf{122}, pp. 1730--1747.

\bibitem[SV]{SV} D.~W.~Stroock, S.~R.~S.~Varadhan (2006). \textit{Multidimensional diffusion processes}. Springer-Verlag, New York.

\bibitem[Sz1]{Sz1} A.~S.~Sznitman (1985). A fluctuation result for nonlinear diffusions. In: \textit{Infinite-dimensional analysis and stochastic processes}. Res. Notes in Math. \textbf{124}, pp. 145--160. Pitman, Boston.  

\bibitem[Sz2]{Sz2} A.~S.~Sznitman (1991). \textit{Topics in propagation of chaos}. \'{E}cole d'\'{E}t\'{e} de Probabilit\'{e}s de Saint-Flour \textbf{XIX}, pp. 165--251. Lecture Notes in Math. \textbf{1464}. Springer-Verlag, Berlin.

\bibitem[Ta]{Ta} H.~Tanaka (1984). Limit theorems for certain diffusion processes with interaction. In: \textit{Stochastic analysis}. North-Holland Math. Library \textbf{32}, pp. 469--488. North-Holland, Amsterdam.

\bibitem[Wa]{Wa} J.~B.~Walsh (1986). \textit{An introduction to stochastic partial differential equations}. \'{E}cole d'\'{E}t\'{e} de Probabilit\'{e}s de Saint-Flour \textbf{XIV}, pp. 265--439. Lecture Notes in Math. \textbf{1180}. Springer-Verlag, Berlin.

\end{thebibliography}
\end{document}